\documentclass[a4paper,11pt]{amsart}
\usepackage{amscd}
\usepackage{calrsfs}

\def\N{\mathbb{N}}
\def\Z{\mathbb{Z}}
\def\Q{\mathbb{Q}}
\def\R{\mathbb{R}}
\def\C{\mathbb{C}}
\def\F{\mathbb{F}}

\def\P{\mathbb{P}}
\def\A{\mathbb{A}}

\def\Fr{\mathrm{Fr}}
\def\tr{\mathrm{trace}}
\def\Ind{\mathrm{Ind}}
\def\sgn{\mathrm{sgn}}
\def\reg{\mathrm{reg}}
\def\H{\mathbb{H}}
\def\Res{\mathrm{Res}}
\def\Spec{\mathop{\mathrm{Spec}}\nolimits}

\def\Sym{\mathop{\mathrm{Sym}}\nolimits}

\def\ker{\mathop{\mathrm{ker}}\nolimits}

\def\H{\mathop{\mathrm{H}}\nolimits}

\def\cycle{\mathrm{cycle}}
\def\cond{\mathrm{cond}}
\def\Art{\mathrm{Art}}
\def\ord{\mathrm{ord}}
\def\GL{\mathrm{GL}}
\def\PGL{\mathrm{PGL}}
\def\SL{\mathrm{SL}}
\def\det{\mathrm{det}}
\def\ad{\mathrm{ad}}
\def\End{\mathrm{End}}
\def\Gal{\mathrm{Gal}}
\def\sw{\mathrm{sw}}
\def\red{\mathrm{red}}
\def\inv{\mathrm{inv}}
\def\unr{\mathrm{unr}}
\def\ss{\mathrm{ss}}
\def\et{\mathrm{\acute{e}t}}

\newtheorem{thm}{Theorem}[section]
\newtheorem{prop}[thm]{Proposition}
\newtheorem{lem}[thm]{Lemma}
\newtheorem{cor}[thm]{Corollary}
\newtheorem{rem}[thm]{Remark}
\newtheorem{dfn}[thm]{Definition}

\newtheorem{exa}[thm]{Example}

\newtheorem{ass}[thm]{Assumption}
\newtheorem{sublem}[thm]{Sublemma}
\newtheorem{claim}[thm]{Claim}

\title {Level optimization in the totally real case}
\author{Kazuhiro Fujiwara}

\begin{document}
\maketitle
\begin{abstract} In this paper, congruences between holomorphic Hilbert modular forms are studied. We show the best possible level optimization result outside $\ell$ for $\ell\geq 3$ by solving the remaining case of Mazur principle when the degree of the totally real field is even.  
\end{abstract}
\section{Introduction}
The aim of this paper is to discuss mod $\ell$ congruences between Hilbert modular
forms. Here $\ell$ is a fixed prime. More precisely, given a holomorphic Hilbert modular
form
$f$ of some level and  weight, we look for another modular form $g$ with a smaller
level and the same weight whose Fourier coefficient are congruent to that of $f$. For
elliptic modular forms, i.e., when $F=\Q$, such studies have been done by J. P. Serre,
B. Mazur, K. A. Ribet, F. Diamond, and by other authors(see the references of \cite{Ri 2}). Amazingly, to
obtain some optimal $g$, the correct condition imposed on $f$ involves an information
from the Galois representation attached to
$f$, so we need to use Galois representations to study congruences. 
\medskip

We
use a representation theoretical terminology in the adelic setting, since it is
essential in the local analysis.  For a
totally real field
$F$, let $I_F =
\{
\iota : F
\hookrightarrow
\R
\}
$ be the set of  the infinite places of $F$. We take an element
$ k = ( k _{\iota} ) _{
\iota
\in I _F}
\in
\Z ^{I_F}$ and  $w \in \Z$, where
$k _{\iota }
\geq 2 $,  and $  k _{\iota} \equiv w  \mod 2 $ for all $\iota$.  Let $
\pi$ be a cuspidal representation of
$\GL _2 (\A_F)
$ having infinity type $(k, w) $ (see
notations for our normalization).  Those types of cuspidal representations are generated
by holomorphic Hilbert cusp forms.

 Fix a prime
$\ell$, and an isomorphism $\C \simeq \bar \Q_\ell$ by the axiom of choice.
 It is known that the finite part $\pi ^{\infty}
$ of $\pi$ is defined over
some algebraic number field, and hence over some finite extension $E_{\lambda}$ of
$\Q_{\ell}$. In this case we say that $\pi$ is defined over $E_\lambda$. For the integer
ring
$\mathcal O _{\lambda}$ of $E_\lambda$, there is a two dimensional continuous
$\lambda$-adic Galois representation 
$$
\rho _{\pi, \lambda} : G _F \to \GL_2 ( \mathcal O _{\lambda})
$$
associated to
$\pi$ (see \cite{O}, \cite{Car 2}, \cite{W}, \cite{T}, \cite{BR} for $F \neq \Q$). 
Here $G_F=
\Gal (
\bar F / F )
$ is the absolute Galois group of $F$.
Two cuspidal representations $\pi$ and $\pi '$ which are
defined over some $E_\lambda$ are {\it congruent} mod $\lambda$ if the
semi-simplifications of the mod
$\lambda$-reductions of their
$\lambda$-adic representations are isomorphic over $k _\lambda$:
$$
(\rho _{\pi,
\lambda} \mod \lambda )^{\ss}
\simeq
(\rho_{\pi ',
\lambda} \mod \lambda  )^{\ss} . 
$$
Here $(-) ^{\ss} $ denotes the semi-simplification. We note that
this notion of congruences is equivalent to the congruences between the Fourier
coefficients of the corresponding normalized Hilbert modular newforms by the
Chebotarev density theorem. 

We call an absolutely irreducible mod $\ell$ continuous representation $\bar \rho : G _F
\to \GL _2 ( \bar k_\lambda )
$ {\it modular} if it is isomorphic to some $\rho _{\pi, \lambda }
\mod
\lambda
$ over $\bar k _\lambda$. Even if $\bar \rho$ is modular, there are many cuspidal
representations giving the same $\bar \rho$, which are all seen as {\it deformations }of
$\bar \rho$ from the viewpoint of Mazur, and the purpose of this paper is restated as to find a
good cuspidal representation $\pi '$ which is optimal for
 modular representation $\bar \rho$ in a suitable sense.  As is already
mentioned,
$F=
\Q$ case is studied well, so we restrict our attention to general totally real
$F$ other than $\Q$.

 Here is our main theorem,
which treats the level optimization at
$ v 
\not \vert \ell $.
\begin{thm}\label{thmA}[Theorem A]
Let
$ \bar \rho : G_F \to \GL _2 (\bar k_\lambda ) $ be a continuous absolutely
irreducible mod
$ \ell$-representation satisfying A1)-A3):
\begin{itemize}
\item [A-1)] $ \ell \geq 3$, and $ \bar \rho \vert _{F (\zeta _{\ell} ) }
$ is absolutely irreducible if $ [F (\zeta _{\ell}) : F]= 2$.

\item [A-2)] $\bar \rho \simeq \rho _{\pi, \lambda} \mod \lambda $. Here  $\pi $ is a
cuspidal representation of $\GL_2(\A _F )$ of infinity type $ (k, w ) $ defined over
$E_\lambda$, satisfying
$ (\pi ^{\infty } ) ^K \neq \{ 0 \} $ for some compact open subgroup $ K =\prod
_u K _u $ of $\GL_2(\A ^{\infty} _F ) $,

\item [A-3)] for a place $ v \not \vert \ell$, $\bar \rho$
is either ramified at $v$, or
$q_v \equiv 1 \mod \ell $.
\medskip

\end{itemize}
Then  there is a cuspidal representation $\pi' $,  having the same
infinity type
$(k, w)$ as $\pi$, defined over a finite extension $E '_{\lambda'}  \supset E_\lambda$
such that the following conditions hold:
\begin{enumerate}
\item The associated $\lambda ' $-adic representation $\rho_{\pi ' , \lambda '}$
gives
$\bar \rho$. 
$ (\pi ^{ ' v,
\infty} ) ^{ K ^v}\neq
\{ 0
\}
$,

\item The conductor 
$\cond (\pi'  _v)
$ of $\pi' _v$ is equal to $ {\Art \bar \rho \vert _{G_v} }$, where $G_v$ is the
decomposition group at $v$, and $\Art
\in
\Z$ means the Artin conductor,

\item  $\det
\rho _{\pi' ,
\lambda  '}\cdot  \chi ^{w +1 } _{\cycle} $ is the Teichm\"uller lift
of
$
\det
\bar
\rho \cdot  \bar \chi ^{w +1} _{\cycle} $.
\end{enumerate}

\end{thm}
If $\pi '$ gives $\bar \rho$, we always have the basic inequality 
$$ \cond \pi _v \geq \Art
\bar
\rho
\vert_{G_v},$$ so the equality in \ref{thmA} (2) is optimal.   When we remove condition A-3), i.e., even
when
$ \bar \rho $ is unramified at $v$ and $ q _v \equiv 1 \mod \ell$, we
still have some 
$\pi '$ giving $\bar \rho$ with the property that
$\pi '_v
$ is spherical, or a special representation twisted by an unramified character.  This missing case in theorem
A, i.e., the case when $\bar \rho $ is unramified at $v$ and $ q_v \equiv 1 \mod \ell $,
was treated by K. A. Ribet in case of $\Q$ \cite{Ri 1}, and by A. Rajaei in the totally
real case.

\medskip
Especially, we note that theorem A includes
\begin{cor}[Corollary A' (Mazur principle)]
 There exists $\pi ' $ as in theorem A when
$\bar
\rho $ is unramified at $v$, and $ q_v \not \equiv 1 \mod \ell$.

\end{cor}
We should note that theorem A is a stronger form of theorem B below, which was obtained
earlier by F. Jarvis
\cite{Ja 1}, \cite{Ja 2} (some additional condition on
$\bar \rho$ is put in the references compared to theorem B, but it is easily
removed, see lemma \ref{lem-auxiliary1} in \ref{subsec-auxiliary}. 
\begin{thm}[Theorem B (Jarvis)] In addition to assumptions A-1)-A-3) of theorem A, 
we assume
\begin{itemize}
 \item [A-4)] If the degree 
$[F:\Q] = g$ is even,  assume that 
$\pi _u $ is essentially square integrable for some $u \neq
v $, $v \not \vert\ell$, with $ K _u = K _1 ( m ^{\cond \pi   _u}_u )$. 
\end{itemize}
Then we have the same conclusion as in theorem A.
\end{thm}
Condition A-4 in theorem B especially excludes the unramified case when the degree
$[F:
\Q ]
$ is even. Logically, theorem A is a consequence of  theorem B and corollary A' in the full form.
Corollary  A' in the even degree case, i.e., when $[F:\Q] $ is even, is one of the main contributions
of this paper. We also give a new proof to theorem B in this article. 
\medskip

Let us explain the method of proofs.  In
\cite{Ja 1}, Jarvis proved a part of corollary A', i.e., the Mazur principle under A-1)-A-4),  by a
detailed study on the arithmetic models of Shimura curves. Note that assumption A-4)
in theorem  B is used to relate $\bar \rho$ to the
cohomology of a Shimura curve associated with a division algebra which is {\it unramified} at
$v$. When $\bar \rho$ is
ramified at $v$
(\cite {Ja 2}), Jarvis used the argument of Carayol in \cite{Car 3}, which does not use any
arithmetic models.
The methods are cohomological.

Our corollary A', especially the Mazur principle in the even degree case, has been thought
difficult since the $\lambda$-adic representation $\rho _{\pi, \lambda }$ may not be obtained
from a Shimura curve in general \cite{T}. Surprisingly, we can use a Shimura curve
attached to a division quaternion algebra which is {\it ramified }at
$v$ in our proof.

Our method to prove corollary A', and to give a new proof to theorem B,  is also
cohomological and summarized in the following way.

For any finite $\mathcal O_\lambda$-algebra $R$, we define a cohomological functor $H_R$ from
Shimura curves, on which
$G_F$ acts, and behaves nicely under any scalar extension $R\to R'$. $H_{\mathcal O _\lambda}
\otimes _{\mathcal O_{\lambda}} E_{\lambda }$ consists of representations $\rho _{\pi ', \lambda '}$'s
giving $\bar
\rho$. The main observation is that the inertia fixed part $H_R ^{I_v}$ also commutes
with scalar extensions, which shows theorem B when $\bar \rho$ is ramified at $v$, since we
have some
$\pi '$ with
$ \rho _{\pi ', \lambda '} ^{I_v} \neq \{ 0\}$ under the assumption $\bar \rho ^{I_v}
\neq \{ 0 \} $. The equality of Artin conductors follows easily from this.

To show the good property of $H _R ^{I_v} $ with respect to scalar extensions in case of
theorem B, we use a deep arithmetic geometrcal result, namely the regularity of the
arithmetic models of Shimura curves using Drinfeld level structures. By this result and the purity theorem of
Zariski-Nagata for
\'etale coverings of a regular scheme, we analyze the cohomology groups directly, without any calculation of
vanishing cycles. Our result is seen as a mod $
\ell$-version of the local invariant cycle theorem.

In the case of corollary A', the analysis of the cohomology is done with the help of a
Cerednik-Drinfeld type theorem obtained by Boutot-Zink \cite{BZ}, based on the principle
used for proving theorem B.
\medskip
Also in \ref{subsec-perfect}, we give an interpretation of Carayol 's lemma \cite{Car 3} by a standard homological
algebra. It is an application of a property of perfect complexes, and also cohomological. 
\medskip
 The methods developed in this article are new even for $F= \Q$. It is our belief that 
a level optimization in easier situation is a consequence of a homological algebra on
Shimura varieties. We hope that the method in this paper is
effective in some higher dimensional cases as well. One may use O. Gabber's absolute
purity theorem (\cite{Fu 1}) instead of the  purity theorem of Zariski-Nagata.

\bigskip

{\bf Acknowledgment}: A part of the work was done during the author's stay at Caltech in
1996, at Universit\'e Paris-Nord in 1998, and at the Institute for Advanced Study for the
academic year 1998-1999.  The author thanks D. Ramakrishnan and J. Tilouine for the
hospitality at Caltech and Paris-Nord.  The author also thanks P. Deligne for the
discussion on some cohomological operations used in the paper. 
 
Finally, the author would like to dedicate this paper to professor Kazuya Kato.  The author has leaned so many mathematical insights from him, and also the beauty of Mathematics. 
\section{Notations} \label{sec-notation}
For a number field $F$, $\mathcal O _F  $ is the integer ring, and $G_F = \Gal (\bar F 
/F)$ is the absolute Galois group. For a place $v$ of $F$, $F_v$ means the
local field at $v$. For a finite place $v$, $o_v$ is the integer ring of $F_v$, with the
maximal ideal
$m_v$,
$k (v ) = o_v/m_v $ is the residue field with the cardinal $q_v$. 
$p_v$ is a uniformizer of $m_v$.
$G_v$ and $I_v$ mean the decomposition and the infertia groups at place $v$,
respectively.

$\A_F$ means the ad\'ele ring of $F$, $\A ^ \infty _F$, and
$(\A _F)_ \infty $ the finite part and the infinite part, respectively. For a
non-zero ideal
$\frak f $ of
$\mathcal O _F
$, we define compact open subgroups of
$\GL _2 (
\A _F ^ {\infty} )$ by  
$$K(\frak f ) = \{ g \in 
\GL _2 (
\prod _{u :\text{finite }}  o _u ),\ g \equiv 1 \mod \frak f  \} ,
$$

$$K _ {11} ( \frak f ) = \{ g \in 
\GL _2 (
\prod _{u :\text{finite }}  o _u ),\ g \equiv
\left(
\begin{array}{cc}
1 &* \\
0&1
\end{array}
\right)
\mod \frak f  \} ,
$$
$$K _ 1 ( \frak f  ) = \{ g \in 
\GL _2 (
\prod _{u  :\text{finite }} o _u ),\ g \equiv 
\left(
\begin{array}{cc}
1 &* \\
0&*
\end{array} 
\right) 
\mod \frak f  \} ,
$$
$$K _ 0 ( \frak f ) = \{ g \in 
\GL _2 (
\prod _{u :\text{finite }}  o _u ),\ g \equiv 
\left(
\begin{array}{cc}
* &* \\
0&*
\end{array}
\right) 
\mod \frak f  \} .
$$
We use the similar notation for the corresponding local group.

For an infinity type $(k, w) $,  $k \in \Z ^I $, $ w \in \Z  $, which satisfies $ k _\iota \equiv w
\mod 2 $, $k ' \in \Z ^I$ is defined by the formula

$$ k + 2 k' =
(w + 2 )
\cdot 
(1, \ldots, 1 ). $$ 

As in the introduction, we fix an isomorphism $ \C \simeq \bar \Q_{\ell}$ by
the Axiom of Choice. 
\medskip
The local Langlands correspondence for $\GL_2 (F_v)$ defines a bijection
between isomorphism classes of $F$-semisimple representation $\rho _v$ and admissible
representation $\pi_v$ of $\GL_2 (F_v)$.

 Our normalization of the local Langlands correspondence is as follows. For the local
class field theory, we assume a geometric Frobenius element corresponds to a
uniformizer. For a finite place $v \not
\vert
\ell$ with spherical
$\pi _v $, a geometric Frobenius element
$
\Fr _ v
$ satisfies
$$ \tr \rho _{\pi , \lambda} (\Fr _v ) = \alpha _v + \beta _v
$$
where $ (  \alpha _v ,\ \beta _v ) $ is the Satake parameter of $\pi _v $ seen as a
semi-simple conjugacy class in the dual group $ \GL_2 ^\vee ( \bar \Q_\ell)$, and
$\pi _v$ is a constituant of the non-unitary induction 
$$
\Ind _{B(F_v)}  ^{G(F_v)} \chi _{\alpha _v, \beta
_v}= \{ f : \GL_2 ( F_v ) \to \bar \Q_{\ell} ,\ f (
\left(
\begin{array}{cc}
a&b\\
0&d
\end{array}
\right)
g ) =\chi _{ \alpha _v, \beta_v } (a, d) \vert a \vert _v f (g)  \} .
$$ Here $B$ is the standard Borel subgroup consisting of the upper triangular matrices, 
$$
\chi
_{\alpha _v,
\beta _v}: B (F_v)
\to   F_v
\times F_v
\to
\bar
\Q_\ell
$$ 
is the unramified character given by $\chi _{\alpha _v, \beta _v}(a, b) = \alpha _v
^{\ord _v a}
\beta _v ^{\ord _v b} 
$.

At infinite places, $\GL_2 (\R)$-representation $D_{ k, w} $ corresponds to the 
{\it unitarily} induced representation
$$\Ind _{B(\R)} ^{G(\R)} ( \mu _{k, w} , \nu _{k, w} ) _{\bold u }
$$
for two characters of split maximal torus
$$\mu _{ k, w} (a) = \vert a \vert ^{\frac  1 2 - k ' } (\sgn a) ^ {-w}
$$
$$
\nu_{ k, w} (d) = \vert d \vert ^{\frac  1 2 -w + k'} 
$$
for $k'$ satisfying $ k -2 + 2 k ' = w $.
\medskip
This normalization, which is the $\vert  \cdot \vert _v  ^{ \frac  1 2 } $-twist of
unitary normalization, preserves the field of definition. The central
character of
$\pi$ corresponds to $
\det
\rho _{\pi,
\lambda} (1)
$. Our normalization is basically the same as in \cite{Car 2}, except one point. In
\cite{Car 2}, an arithmetic Frobenius element corresponds to a uniformizer.

The global correspondence $\pi \mapsto \rho _{\pi , \lambda} $ is compatible
with the local Langlands correspondence for $ v \not \vert \ell$ (\cite{Car 2}
th\'eor\`eme (A), see \cite{W}, \cite{T}, theorem 2, for the missing even degree cases):
 the 
$F$-semisimplification of
$
\rho_{\pi , \lambda} \vert _{G_v}
$ corresponds to $\pi_v$ by the local correspondence normalized as above.

\section{Preliminaries}\label{sec-preliminary}
\subsection {Shimura curves and Hida varieties} \label{subsec-shimura}

We assume $[ F : \Q ] > 1 $ in the
following (cf. introduction). One can include
$\Q$-case with a slight modification. Fix an element 
$\iota _1
\in I _F
$. Take a division quaternion algebra $D$ over $F$ which ramifies at all
infinite places other than $ \iota _1 $, and ramifies possibly at $\iota_1$. We fix an
identification
$$
D \otimes _{\Q }\R \simeq M _2 (\R ) ^ { g '} \times \H ^
{ g - g '} ,
$$
where $ g ' = 1 $ or $0$ according to $ D$ is split at $\iota _ 1$ or not. Here $\H$ is
the Hamilton quaternion algebra over $\R$. We put 
$$
G_D = \Res _{ F/ \Q} D ^\times.
$$
Here $\Res$ means the Weil restriction of scalars. 
For a compact open subgroup $ K \subset D ^ {\times} ( \A ^ {\infty
} ) $, we define the associated modular variety with a complex structure, by
$$
S_ { K } = D ^ {\times }\backslash D ^ {\times} ( \A ) / K
\times  K _ {\infty } .
$$
Here $K _{\infty} $ is the maximal compact subgroup of $ D ^ {\times} ( \R ) $
modulo center. When $g'$ is one, $S_K$ is a Shimura curve, and it has a canonical
model $S_{K, F} $ defined over $F$. When $g'$ is zero, the zero-dimensional variety was
used extensively by H. Hida \cite{H 1} in his study of $\ell$-adic Hecke algebras. Note that this
is {\it not} a Shimura variety in the sense of Deligne. 
\medskip
\subsection{ $\lambda$-adic local systems} \label{subsec-local}
Fix a prime $\ell$, and an
$\ell$-adic field $E_\lambda$ with the integer ring
$\mathcal O _\lambda$. We denote the maximal ideal of $\mathcal O_\lambda$ by $\lambda$ and the
residue field by
$k_\lambda$. 

For a pair $(k, w)$, $ k \in \Z ^{I_F}$, $w \in \Z$ as in the introduction, we define an
$\ell$-adic sheaf
$\bar
{\mathcal {F}} _{k,w}
$ on $S_K $ (\cite{Car 2}, p.418-419 for $g'=1$). 

Since we need a $\Z_{\ell}$-structure for this sheaf, we briefly review the construction, assuming $D$
splits at all $v \vert \ell$, which is sufficient for our later use.

Let $\pi _{\ell} : \tilde S _{\ell} \to S _{D, K} $ be
the Galois covering corresponding to $\prod _{u \vert \ell} K _u / K \cap \overline {F
^{\times}} $. Here $  \overline {F^{\times}}$ is the closure inside $D ^{\times} ( \A ^f _F) $.

We choose an
isomorphism 
$ D \otimes _\Q \Q_{\ell} \simeq \prod _{v \vert \ell} M_2 (F_v)$. This determines a
maximal hyperspecial subgroup of $D ^ \times (\Q_\ell)$, which we identify with $\prod
_{v \vert \ell }\GL_2 (o_v)$. We take
$ E_\lambda$ large enough so that all
$F_ v ,\ v\vert\ell$ are embeded into
$E_\lambda$ over
$\Q_\ell$. So the representation 
$$ 
V_{k, w} = \otimes
_{\iota \in I _F} (\iota \det) ^{-k' _{\iota}}\Sym ^{  \vee \otimes( k _{\iota}-2)} 
$$ 
of $ D ^{\times}
(\Q _{\ell})$ is defined over $E_\lambda$, and has an
$\mathcal O_\lambda$-lattice 
$$
V_{ ( k,w ) , \mathcal O_\lambda} = 
\otimes _{v \vert \ell, \iota :F_v \to E_\lambda }(\iota \det) ^{-k' _{\iota}}\Sym ^{ 
\vee \otimes (k _{\iota}-2)}
\mathcal O _v  ^{\oplus 2}) \otimes _{\iota} \mathcal O _\lambda
.$$
 The $\mathcal O _\lambda$-smooth sheaf $ \bar {\mathcal { F} } ^D _{k, w}$
is obtained from the covering $\pi _{\ell}$ and the representation $V_{ ( k,w ) , \mathcal
O_\lambda}$ by contraction. Note that the action of $ K \cap
\overline{ F^\times }$ on the representation  $V_{k, w} $ is trivial if $K$ is
sufficiently small, so that the sheaf is well-defined. 

See
\cite{Car 2}, p.418-419 for the Betti-version of $\bar {\mathcal {F}} _{k, w} $ (note that we
have taken the dual of the sheaf in the reference). By the comparison theorem in
\'etale cohomology, those two cohomologies are canonically isomorphic, so we do not make
any distinction unless otherwise specified. 
\medskip
The projective limit  $S = \varprojlim _K S_K $ admits a right $D ^\times (\A
^\infty )
$-action by the multiplication from the right, and 
the $E_\lambda
$-sheaf 
$$
 \bar { \mathcal {F}} ^D _{k, w,\ E_\lambda }=  \bar {\mathcal {F}} ^D _{k, w}\otimes _{\mathcal O_\lambda
} E_\lambda 
$$
is $ D ^\times (\A ^\infty ) $-equivariant by the
construction.

The lattice structure  $\bar{ \mathcal {F} }^D _{k, w }$is preserved by the $D
^\times(\A ^{\ell , \infty}  )  $-action.

When $g' $ is one, the sheaf $ \bar { \mathcal {F}} ^D _{k, w}$ is canonically defined over $F$ by
the theory of canonical models, which we denote by $ \mathcal F ^D _{k, w}$. By \cite{Car 1}
2.6,
$
\mathcal F ^D _{k, w}$ is pure of weight $w$. 
This canonical
$F$-structure gives a continuous $G_F$-action on the cohomology, giving the
decomposition
$$
H ^1 _{\et}( S_{ K, \bar F} ,\ \bar{ \mathcal {F}}_{ k, w}  \otimes E_{\lambda }) = \oplus
_{\pi }
\rho _{
\pi,
\lambda }\otimes _{\mathcal O _\lambda }  (\pi ^{\infty} ) ^{\tilde K } 
$$
for cuspidal representations $\pi $ of $D ^\times (\A_F) $ with infinity type $ ( k, w)
$ which does not factor through the reduced norm (\cite{Car 2}, 2.2) assuming that
$E_\lambda
$ is sufficiently large. By infinity type
$ (k, w)
$, we mean that $ \pi_\infty$ has the form
$$\pi _\infty =D_{ k _{\iota _1}, w}
\otimes _{ v \in I_F, v \neq \iota _1} \bar D_{ k _v, w} $$ (see \cite{Car 2} \S0 for
the notation $\bar D_{ k_v, w} = (\iota \det) ^{-k' _{\iota}}\Sym ^{  \vee \otimes k
_{\iota}-2}$, which corresponds to $ D_{ k_v, w} $ by the Jacquet-Langlands
correspondence \cite{JL}).

 \subsection{Hecke algebras and correspondences}\label{subsec-Hecke}
We keep the same notations as in \ref{subsec-local}. We assume that $K$ is locally factorizable, i.e., $ K =\prod
_u K_u
$, and
$D$ is split at all
$ v
\vert \ell$. For a finite set of finite places $\Sigma $ containing  $ \{ v; \ v \vert \ell \} $,  let
$H ( D ^{\times} (\A ^{\Sigma ,
\infty}) ,\ K ^{\Sigma }) $ be the
convolution algebra  over $\mathcal O _\lambda $. Namely,  $H ( D
^{\times} (\A ^{\Sigma, 
\infty}) ,\ K ^{\Sigma }) $ is the set of the compactly supported
$\mathcal O _\lambda
$-valued smooth
$K ^{\Sigma }$-biinvariant functions on $D ^\times (\A ^{\Sigma, \infty}) $. The algebra
structure is given by the convolution.

We consider the action of $H ( D ^{\times} (\A ^{\ell,
\infty}) ,\ K ^{\ell}) $ on the cohomology of
$S_{\tilde K }$ induced from the adelic right action on $S= \varprojlim _{\tilde K } S _{\tilde
K } $. We briefly review the basic facts since the relationship with the Verdier duality is subtle
and important for our purpose.

For two compact open subgroups $K ,\ K '   \subset D ^\times(\A^\infty) $, $g \in
D ^\times(\A^\infty)$ define an algebraic correspondence
$[KgK '] 
$: the first projection
$S_{ K
\cap g ^{-1} K ' g }
\to S_K
$ and the second 
$S_{ K
\cap g ^{-1} K ' g } \to S_{ g ^{-1} K ' g } \simeq  S_{K'}$. 

The 
correspondence induced by $K g K '$ from $S_K $ to $S_{K'} $ is dual to 
$ K ' g ^{-1} K $ from $S_{K'} $ to  $S_K $ by the definition.

Since $\mathcal F_{k, w, E_\lambda } $ is $D ^{\times} (\A ^\infty)$-equivariant by the
construction,  $[KgK'] $ gives
$$
[KgK']: R\Gamma ( S_{ K } ,\ \bar {\mathcal F} _{k, w} ) \otimes _{\mathcal O _\lambda }E
_\lambda  \to R\Gamma ( S_{K '} ,\
\bar
{\mathcal F } _{k, w} ) \otimes _{\mathcal O _\lambda }E
_\lambda  $$ 
via cohomological correspondences they define.

We fix a uniformizer $p_v$ of $F_v$ at each finite place $v$. If the $v$-component
$g_v$ of
$ g
\in D ^\times (\A ^\infty ) $ satisfies

$$
 \GL_2 (
o_v )  g _v
\GL_2 ( o_v )\text{ is represented by }
\left(
\begin{array}{cc}
 p _v ^a &0 \\
0&p _v ^ b
\end{array}
\right)
\text{ with } a,\ b \geq 0  \text{ for all }v \vert \ell 
$$
the $\mathcal O _\lambda$-lattice structure $\mathcal F _{k, w}  $ is preserved, and 
$$
[KgK']: R\Gamma ( S_{ K } ,\ \bar{ \mathcal {F}} _{k, w} )  \to R\Gamma ( S_{K '} ,\
\bar{\mathcal {F}} _{k, w} ) 
$$ 
is induced. We call this the {\it standard }action of $[KgK'] $. Especially, we have the action of
$H ( D ^{\times} (\A ^{\Sigma,
\infty}) ,\ K ^{\Sigma}) $ for $ K = K _\Sigma \cdot K ^\Sigma $ by extending $ K ^\Sigma g
K ^\Sigma$ to
$ K \tilde g K $, $\tilde g = ( (1 _{D(F_v)} ) _{ v \in \Sigma} ,  g ) ,\ g \in D ^\times
(\A ^{\Sigma, \infty})
$.
\medskip

It is a natural question to ask whether there is a complex  of
$H ( D ^{\times} (\A ^{\Sigma,
\infty}) ,\ K ^{\Sigma}) $ -modules which represents $R\Gamma ( S_{ K } , \bar {\mathcal {F}} _{k,
w} ) 
$ or not. The answer is affirmative, which we state it in the form of proposition. In section \ref{sec-proofB}, it becomes quite important.
\begin{prop}\label{prop-Hecke1}
There is a complex $L ^\cdot  $ of $ H ( D ^{\times} (\A ^{\Sigma,
\infty}) ,\ K ^{\Sigma}) $-modules bounded from below which represents $R\Gamma ( S_{
K } , \bar {\mathcal {F}} _{k, w} )  $ in $ D ^+ ( S_{
K } , \mathcal O _\lambda ) $. The induced action $ of  H ( D ^{\times} (\A ^{\Sigma,
\infty}) ,\ K ^{\Sigma })  $on $ H ^ q ( S_{
K } , \bar { \mathcal {F} }_{k, w} )   $ coincides with the induced action by the Hecke
correspondences.  
\end{prop}
\begin{proof} We work with the Betti version for simplicity (this case is sufficient for our
later purpose). Let $ L^\cdot
$ be the Godement's canonical resolution of $\bar{  \mathcal {F}}_{ k, w} $. Since all the maps
in defining the algebraic correspondences are finite and \'etale, cohomological
operations defining the action of a cohomological correspondence are actually defined on
$L^\cdot $: for example, for a finite \'etale morphism $ f : X \to Y$, the trace map 
$ f _ ! f ^! \mathcal F  = f _* f ^* \mathcal F \to \mathcal F  $ for  a sheaf $\mathcal F$ is defined by
$$
\sum _{ y \in f ^{-1} (x) } \mathcal F _y \to \mathcal F _x ,\ ( f_y ) _{ y \in f ^{-1} (x)
} \to \sum _{ y \in  f ^{-1} (x) } f_y 
$$  
on the level of stalk, which gives a lift of the trace map to Godement resolution
$L^\cdot$.  In this way one has an action of the convolution algebra on
$\Gamma (L^\cdot ) 
$, and $ R \Gamma (S_K , \ \mathcal F_{k, w} ) $ belongs to the derived category of $ H (D
^\times (\A ^{\Sigma, \infty }) ,\ K  )
$-modules bounded below. Especially, for
$K =
\prod K_v$, the actions at two different places commute. 

\end{proof} 
\begin {rem}\label{rem-Hecke1}
For the \'etale cohomology with finite coefficients, the claim is proved by the same argument. 
\end{rem}

\subsection{Duality formalism}\label{subsec-duality}
Take a finite set of finite places $\Sigma $ containing $ \{v;\  v \vert \ell\} $.
Assume
$
\tilde K =
\prod _{u }
\tilde K _ u
$, and
$D$ is split at all finite $ v \not \in \Sigma $. 

Note that 
$$
\bar {\mathcal F} ^{\vee}  _{k, w}\simeq \bar{ \mathcal {F}} _{ k, -w}. \leqno{(\dagger)}
$$
 giving a perfect pairing in the derived category of $\mathcal O_\lambda$-modules
$$ 
R
\Gamma  (S_{K } ,\
\bar
{\mathcal F}   _{k, w} )\otimes ^{\mathbb L} _{\mathcal O_\lambda} R
\Gamma  (S_{K }  ,\ \bar{ \mathcal F}   _{k, -w}) \to \mathcal O _\lambda ( -g' )[-2g']  
$$ 
by the Poincar\'e duality.
For the 
$D^\times(\A ^{\Sigma, \infty} )$-action, if we consider the standard action of $ g$ on $
 R
\Gamma  (S_{K } ,\
\bar
{\mathcal F}   _{k, w} )$, by the Poincar\'e duality, this corresponds to the standard action of $g
^{-1} $
on $ R
\Gamma  (S_{K }  ,\ \bar {\mathcal {F} }  _{k, -w}) $ since the isomorphism $(\dagger)$ sends $g$-action
to $ g ^{-1}$.  
\bigskip

 For the relation between $H ( D ^{\times} (\A ^{\Sigma, \infty}) ,\
 K ) $-action and the Verdier duality,
we have the following proposition by discussion of section \ref{subsec-Hecke}. 

\begin{prop}\label{prop-duality1} The standard action $ R \Gamma  (S_K ,\
\bar {\mathcal {F} }_{k, w}) \to R \Gamma  (S_{K'} ,\
\bar{ \mathcal {F}} _{k, w}) $ induced by $ [K g K ']$ is dual to the standard action $ R \Gamma 
(S_{K'} ,\
\bar { \mathcal {F }}_{k, -w})(g') [2g'] \to R \Gamma  (S_{K} ,\
\bar {\mathcal {F}} _{k, -w})(g')  [2g']$ by $ [K
' g ^{-1} K ]$. 
\end{prop}
\medskip
We have two geometric actions of the convolution algebra, which we
call the {\it standard} action and the {\it dual }action. By the dual action of $ [ K g K ]$ on
$  R
\Gamma  (S_K ,
\bar {\mathcal { F}} _{k, w}) $, we mean the standard action of $ [ K g^{-1} K ]$.

By proposition \ref{prop-duality1}, the standard action of $T_\Sigma$ becomes the dual action
by the Poincar\'e duality.  
\bigskip

We define standard Hecke operators. We choose a unformizer $p_v$ of $F_v$ for any finite place
$v$. 
$a (p_v)$ (resp. $ b ( p_v )$) is the element in $ D ^\times ( \A ^\infty )$ having 
$
\left(
\begin{array}{cc}
1 & 0\\ 
0 & p _v \\
\end{array}
\right)
$ (resp. 
$ 
\left(
\begin{array}{cc}
p_v & 0\\ 
0 & p _v \\
\end{array}
\right)
$) as it's $v$-component,
and 1 as other components.

Then the Hecke algebra $ T _\Sigma= H ( D
^\times ( \A ^{\Sigma , \infty }) ,\ K ^\Sigma ) 
$
over $\mathcal O _\lambda$ is isomorphic to the $\mathcal O _\lambda$-algebra generated by
indeterminates
$ [ T _v]$,
$[T_{v, v}],\ v \not \in \Sigma
$, adding $ [T_{v, v}]^{-1} $ for $ v \not \in \Sigma $. 
Here $[T_v]  $ is given by the characteristic function of  $
\tilde K ^\Sigma a( p_v) 
\tilde K^\Sigma $, $[T_{v, v} ] $ is given by  the characteristic function of $ \tilde K ^\Sigma b
(p_v) 
\tilde K ^\Sigma $.
For a $T_\Sigma $-module $M$, we define the dual $T_\Sigma$-action on $M$ by
$$
\begin{cases} [T_v] \mapsto [T_{v, v} ] ^{-1} \cdot  [T_v]\\
[T_{v, v}] \mapsto [ T_{v, v} ] ^{-1} .
\end{cases} 
$$
\begin{exa}\label{exa-duality1} For any continuous mod $\ell$ Galois representation $\bar  \rho :
G_\Sigma
\to
\GL_2 (\bar k_\lambda )$, we define a maximal ideal $m_{\bar \rho}$ of $T_\Sigma$ by 
$$
[ T _v] \mapsto \tr \bar \rho (\Fr_v),\ [ T_{v, v}] \mapsto q_v ^{-1} \det \bar \rho
(\Fr_v) .
$$
The maximal ideal corresponding to the dual action is $m _{ (\det \bar  \rho) ^{-1} 
\otimes
\bar \rho (1)}
$.
\end{exa}

In case of Shimura curves, the definition is compatible with the $T_\Sigma$-action on
cohomology groups:
 
By choosing the canonical resolution, $ R \Gamma (  S_{\tilde K } , \ \bar {\mathcal
F} _{k, w} )$ belongs to the derived category of $T_\Sigma$-complexes bounded below,
sending $[T_v] $ and $[T_{v,v}] $ to the standard actions of $[ \tilde K 
\left(
\begin{array}{cc}
1 & 0\\ 
0 & p_v \\
\end{array}
\right)
\tilde K ]$ and
$[ \tilde K
\left(
\begin{array}{cc}
p_v & 0\\ 
0 & p_v \\
\end{array}
\right)
\tilde K ]$. 
Then the dual action defined by the Poincar\'e duality is given by the above dual
$T_\Sigma$-action since they are defined by $[ \tilde K
\left(
\begin{array}{cc}
1 & 0\\ 
0 & p^{-1} _v \\
\end{array}
\right)
\tilde K ]$ and
$[ \tilde K
\left(
\begin{array}{cc}
 p _v^{-1} & 0\\ 
0 &  p _v^{-1} \\
\end{array}
\right)
\tilde K ]$.  
\bigskip

\subsection {Modules of type $\omega$}\label{subsec-type}
We keep the notations, especially  $T_\Sigma =  H (  D ^\times ( \A ^{\Sigma , \infty } ),\
K^\Sigma )
$. We define a class of
$T_\Sigma$-modules. 
\begin{dfn}\label{dfn-type1} Consider the category $\mathcal C _{T_\Sigma}$ of
$T_\Sigma$-modules which are finitely generated as $\mathcal O_\lambda$-modules. 
We call an object $N $ in $\mathcal C _{T_\Sigma}$ of type $\omega$
\begin{enumerate}
\item if it has a finite length, then for any constituant $N'$ appearing in
the Jordan-H\"older sequence 
$$ 
[T_v ]^2 = [T_{v,v} ] (1 + q_v)^2
$$ 
holds on $N '$ and 
for
$v\not \in \Sigma$,
\item in general $N$ is of type $\omega$ if and only if the graded modules $ \lambda
^n N/
\lambda ^{n+1} N$ ($n \in \N$) for the $\lambda$-adic filtration are all of type
$\omega$.

\item A maximal ideal $m$ of $T_\Sigma$ is of type $\omega$ if $T_\Sigma /m$ is of
type $\omega$.
\end{enumerate}

 By
$\mathcal C _{\Omega}
$ we denote the subcategory of $ \mathcal C_{T_\Sigma} $ consisting of the $T_\Sigma $-modules
of type
$\omega$.
\end{dfn}
It is easy to see that $\mathcal C _\Omega$ forms a Serre subcategory of $\mathcal C _{T_\Sigma} $, and is
stable under the dual action of $T_\Sigma$. By
$
\mathcal C _{N\Omega }$, we mean the quotient category of $\mathcal C _{T_\Sigma}$ by $\mathcal C _\Omega
$.\par
\medskip
A typical example of modules of type $\omega$ is obtained by a one dimensional
representation $\chi : D ^\times ( \A ^{\Sigma , \infty } ) \to E ^\times _{\lambda}  $.
The induced $ T_\Sigma$-action gives
a module of type $\omega$.  

\bigskip
The following proposition shows that the maximal ideals of $T_\Sigma$ of type 
$\omega$ correspond to very special reducible representations. 

\begin{prop}\label{prop-type1}
Assume that a continuous representation $\bar  \rho : G_\Sigma \to \GL_2 (\bar k_\lambda
)$ satisfies 
$$\tr \bar \rho (\Fr_v)
^2 = ( 1+ q_v ) ^ 2  q_v ^{-1} \det \bar \rho (\Fr_v) \leqno{(\dagger)}
$$
 for all $v \not \in \Sigma$. Then $\bar \rho$ is reducible, and the semi-simplification
$\bar \rho ^{\ss}$  satisfies 
$$
\bar \rho ^{\ss} \simeq \bar \chi \oplus \bar
\chi (1)
$$ 
for some one dimensional character $\bar \chi: G_\Sigma \to k_\lambda ^\times$. In other
words,
the maximal ideal $m_{\bar \rho}$ corresponding to $\bar \rho$is of type
$\omega$ if and only if $\bar \rho ^{\ss} \simeq
\bar \chi \oplus \bar \chi (1)  $ for some $\bar \chi$.
\end{prop}

\begin{proof} Let $ \omega =
\chi_{\cycle} \mod \ell$ be the Teichm\"uller character.  By the Chebotarev density
theorem, \ref{prop-type1} $(\dagger)$ is equivalent to 
$$
\tr (g ;  \ad ^0 \bar \rho ) = \tr (g; 1 \oplus \omega \oplus \omega ^{-1} ) 
$$
for any $ g \in G_\Sigma$. 
Since $ \ad ^0 \bar \rho$ satisfies $ \Lambda ^3 \ad ^0 \bar \rho\simeq k $ and
self-dual, the above equality of the traces implies the equality of characteristic polynomials
$$
\det (1-gT ;  \ad ^0 \bar \rho  ) = \det (1-gT ; 1 \oplus \omega \oplus \omega ^{-1} ) 
$$
for any $ g \in G_\Sigma$. 
By the Brauer-Nesbit theorem, the equality of semi-simplifications 
$$
(\ad  ^0 \bar \rho ) ^{\ss} \simeq 1
\oplus \omega \oplus \omega ^{-1}
$$ 
as $G_\Sigma$-modules follows. If $\bar \rho$ is irreducible, the only possibility is
to have
$\omega $ or $\omega ^{-1} $ as a submodule of $\ad ^0 \rho $, and hence
$[F(\zeta_\ell) : F ] = 2 $ and $\rho$ is induced from $F(\zeta_\ell) $. But this type of 
representations do not have the adjoint representation of the above form. 
We conclude that $\bar \rho$ is reducible, and the claim follows.  

\end{proof}

\medskip
\begin{rem} \label{rem-type1} The notion of modules of type $\omega$ is stronger than the
notion of Eisenstein modules in
\cite{DT 1}.  
\end{rem}
\medskip
\subsection{Cohomological lemmas}\label{subsec-cohomology}
Let $S$ be a strict trait with the generic point $\eta$ and the closed point $s$.
$I = \Gal ( \bar \eta / \eta)$. Here $\bar  \eta $ is a geometric point over $\eta$. 

\begin{lem} \label{lem-cohomology1} Let $ X$ be a flat regular scheme of finite type over $
S$.  Assume $\ell$ is prime to the residual characteristic of $S$, and take a
coefficient ring $\Lambda $ which is finite over $\Z _{\ell} $. For a $\Lambda
$-smooth sheaf
$\mathcal F$ on $X$, there are two exact sequences 
\begin{itemize}
\item [(a)] $$ 0 \to H ^0_{\et} ( X_{\bar \eta} , \ \bar {\mathcal {F}} ) ( -1) _I \to H ^1_{\et}( X_{\eta} ,\
\mathcal F _{\eta}) \to H ^1 _{\et}( X_{\bar \eta}  ,\ \bar {\mathcal F} ) ^I \to 0 ,  $$

\item[(b)] 
$$ 
 0 \to H ^1_{\et}( X , \ \mathcal F) \to H ^1 _{\et}( X_{\eta} ,\ \mathcal F _{\eta}) \to \prod _{Y \in
J }H ^0 _{\et}(Y _{\reg},\
\mathcal F \vert _{Y _{\reg}})  (-1) .
$$ 
\end{itemize}
Here $J$ is the set of irreducible components of $X _s $, and $(-)_{\reg} $ means
the regular locus of $(-)$.
\end{lem}

\begin{proof} 
The exactness of  \ref{lem-cohomology1} (a)  follows from the Hochshield-Serre spectral sequence applied to
the morphism from $X_{\bar
\eta}$ to  $X_\eta
$. We show  \ref{lem-cohomology1} (b).  $H ^0_{\et}
(Y _{\reg},\
\mathcal F \vert _{Y _{\reg}})  (-1) $ is canonically isomorphic to $ H ^1 _{\et}(
\operatorname {Frac} o _{X, \bar \eta _Y } ,\ \mathcal F _{\bar \eta_Y} )$, where $\bar
\eta_Y$ is a geometric generic point of
$Y$. Then the exactness of  \ref{lem-cohomology1} (b) is equivalent to the following claim.
\begin{claim}\label{claim-cohomology1}
An $ \mathcal F_{\eta} $-torsor over $X_{\eta}$ which becomes unramified at all maximal
points of
$X_s$ extends to an $\mathcal F$-torsor over $X$ uniquely up to isomorphism. 
\end{claim}
\ref{claim-cohomology1}
is a consequence of the purity theorem of Zariski-Nagata \cite{G}, which says
that, if
$X$ is regular, the category of the \'etale coverings of $X$ is equivalent to the
subcategory of the
\'etale coverings of $X_{\eta} $ which becomes unramified at all maximal points of $X_s$.
\end{proof}

\medskip
Fix a commutative $\mathcal O_\lambda $-algebra $A$. We assume the
following conditions on $(A, X, \mathcal F)$.
\begin{ass} \label{ass-cohomology1} For any $ i \in \Z$, 
\begin{enumerate}
\item $H ^i _R = H ^i _{\et}( X_{\bar
\eta} ,\
\bar
{\mathcal F}
 \otimes _{\mathcal O_\lambda } R )$, $ \mathcal H ^i _R = H ^i _{\et}( X,\
\mathcal F 
 \otimes _{\mathcal O_\lambda } R )$ are $(A, R)$-bimodules for any commutative finite
$\mathcal O_\lambda$-algebra $R$ which commutes with $I$-actions.
\
\medskip
\item $H ^i _R \to H ^i _{R'} $ is an
$A$-module homomorphism for any $\mathcal O_\lambda$-algebra homomorphism $ R \to R' $.
\medskip
\item All homomorphisms in the long cohomology exact sequence $\{ H ^i _R \} _{i \in \Z}
$ are
$A$-module homomorphisms, functorial for  any ring extension $R \to R'$.
\end{enumerate}
\end{ass}

\begin{lem}\label{lem-cohomology2}
Fix a maximal ideal $m$ of $A$. Let $\mathcal E$ be a Serre subcategory of the category of $A$-modules $(A\text{-mod})$,
consisting of 
$\mathcal O_\lambda$-modules of finite length satisfying the following conditions. 
\begin{enumerate}
\item $\mathcal E$ is stable under
$(-)
\otimes _{\mathcal O_\lambda } R$ for any finite
$\mathcal O_\lambda$-algebra $R$.
\item The localization $M_m$ at $m$ of any
element $M$ in $\mathcal E$ vanishes.

\item $ X$ is a proper flat regular curve over $ S$ with
the smooth generic fiber, and $\mathcal F$ is an $\mathcal O_\lambda$-smooth sheaf on $X$.
Assumption \ref{ass-cohomology1} is satisfied for $(A, X, \mathcal F) $.
\item $ H ^0 _{\et}(X _{\bar \eta} ,\ \bar{\mathcal F} \otimes _{\mathcal O_\lambda} R ),\ H
^2_{\et} ( X_{\bar
\eta},\
\bar
{\mathcal F}  \otimes _{\mathcal O_\lambda }R)
\in
\mathcal E
$ for any $\mathcal O _\lambda$-module
$R$ of finite length.\par
\item $H ^0_{\et} ( X_s,\ \mathcal F \vert _{X_s} \otimes _{\mathcal O_\lambda} R) ,\
H ^2_c ( Y_{\reg},\ \mathcal F \vert _{Y_{\reg}} \otimes _{\mathcal O_\lambda} R ) \in
\mathcal E$ for any irreducible component
$Y$ of $X_s$ and for any $R$ of finite length. 
\end{enumerate}
Put $ H  _R = H ^1 _R$, $\mathcal H _R = \mathcal H^1  _ R$.
Then the following hold:
\begin{itemize}
\item[(a)] $H_{R, m} $ is $\mathcal O_\lambda $-free, and  $ H  _{R , m}
\otimes  _{\mathcal O_\lambda} k_\lambda
\simeq H _{R/ \lambda R , m}
$ if $R$ is $\mathcal O_\lambda$-flat. 

\item[(b)] $\mathcal H_{R, m} $ is $\mathcal O_\lambda $-free, and $ (\mathcal H  _{R, m} )
\otimes _R k_\lambda 
\simeq \mathcal H  _{R/ \lambda R , m}$ if $R$ is $\mathcal O_\lambda$-flat.

\item [(c)]$ \mathcal H _{R, m} \simeq (H
_{R, m}) ^I
$ if $R$ is $\mathcal O_\lambda$-flat.  
\end{itemize}
\end{lem}
\begin{proof} First note that $M _m $ is finitely generated as an
$\mathcal O_\lambda$-module for any $A$-module $M$, finite over $\mathcal O_\lambda$, since the $A$-action factors
through a subalgebra in $\End _{\mathcal O_\lambda } M
$, and hence a finite $\mathcal O_\lambda$-algebra, which is a product of local rings. 

We prove (a). By
the exact sequence
$$
H ^0_{\et} ( X_{\bar \eta } ,\ \bar {\mathcal {F}} \otimes _{\mathcal O_\lambda} R/ \lambda R ) \to H
^1 _{\et}(  X_{\bar \eta } ,\ \bar { \mathcal {F}} \otimes _{\mathcal O_\lambda} R ) \overset \lambda \to
 H ^1 _{\et}(  X_{\bar \eta } ,\ \bar {\mathcal F}  \otimes _{\mathcal O_\lambda} R ) 
$$
$H _{R , m}$ is an $\mathcal O_\lambda $-flat module since $H ^0_{\et} ( X_{\bar \eta } ,\
\bar{ \mathcal F }\otimes _{\mathcal O_\lambda} R/ \lambda R )$ vanishes after localization at $m$. 
$$
 H ^1 _{\et}(  X_{\bar \eta
} ,\ \bar {\mathcal F} \otimes _{\mathcal O_\lambda} R ) \overset \lambda \to H ^1
_{\et}(  X_{\bar \eta } ,\ \bar {\mathcal F} \otimes _{\mathcal O_\lambda} R ) \to  H ^1_{\et} ( 
X_{\bar
\eta } ,\ \bar {\mathcal F } \otimes _{\mathcal O_\lambda} R/ \lambda R ) \to H ^2 _{\et}( X_{\bar \eta
} ,\
\bar {\mathcal F} \otimes _{\mathcal O_\lambda} R/
\lambda R )
$$
is exact, and the claim for $H _R $
follows.

For (b), the $\mathcal O_\lambda$-freeness is proved similary as above once we know that  $H
^0_{\et} ( X_s,\ \mathcal F \vert _{X_s} \otimes _{\mathcal O_\lambda} R) \in \mathcal E $. This follows
from 
$$
H ^ 0 _{\et}( X_s, \mathcal F \vert _{X_s} \otimes _{\mathcal O_\lambda} R) ) \overset {(1)} = H ^ 0
_{\et}(X,
\mathcal F 
\otimes _{\mathcal O_\lambda} R)  \subset H ^ 0 _{\et}(X_\eta ,\mathcal F _\eta  \otimes _{\mathcal
O_\lambda} R ) 
$$
and by (2). Here $(1)$ is an isomorphism by the proper  base  change theorem.
 For the rest of (b), it suffices to see $  H ^2 _{\et}( X_s , \
\mathcal F \vert _{X_s}\otimes _{\mathcal O_\lambda} R/ \lambda R ) \in \mathcal E$.
This follows from 
$$ H ^2 _{\et}( X_s , \
\mathcal F \vert _{X_s}\otimes _{\mathcal O_\lambda} R/ \lambda R ) =  
\oplus _{Y \in J} H ^2 _c ( Y_{\reg},\
\mathcal F \vert _{Y_{\reg} }\otimes _{\mathcal O_\lambda} R/ \lambda R)  .
$$

For (c), consider the commutative diagram
$$\CD
 \mathcal H _{R , m } \otimes _{\mathcal O_\lambda} k_\lambda @>(1)>> H^I _{R, m}\otimes _{\mathcal
O_\lambda } k_\lambda \\ @V (2) VV @VV (4)V \\
\mathcal H _{R/\lambda R , m} @>(3)>> H ^I _{R / \lambda R , m} .\\
\endCD
$$
  (2) is an isomorphism by (b). (3) is an isomorphism by
lemma \ref{lem-cohomology1}.  The composition 
$$ 
H^I _{R, m}\otimes _{\mathcal O_\lambda } k _\lambda \overset {(4)} \to H ^I _{R
/
\lambda R , m}
\subset H  _{R / \lambda R , m} = H _{R, m} \otimes _{\mathcal O_\lambda } k_\lambda 
$$ 
is
injective since
$H ^I _{R, m } \subset H _{R, m}$ is an $\mathcal O_\lambda$-direct summand by the definition
and (a). So the map
$(1)$ is an isomorphism. $ \mathcal H _{R , m} \simeq H^I _{R, m} $ follows since
both modules are $\mathcal O_\lambda$-finite free. 
\end{proof}
\begin{rem}\label{lem-cohomology1}
Since $\mathcal C$ is a Serre subcategory, it suffices to check the conditions for  $ R = A/m  $. 
\end{rem}

\section {Proof of theorem B in the ramified case}\label{sec-proofB}
\subsection{Auxiliary places} \label{subsec-auxiliary}
Technically, we need to choose some auxiliary place $y$ such that the discrete subgroups are torsion-free, and the
Hecke algebra does not introduce essentially new component at $y$, the idea introduced by
F. Diamond and R. Taylor (\cite{DT 2}, lemma 11). The extra assumption on $\bar \rho$
when
$[F(\zeta _{\ell}):F] =2
$ is necessary to make this change possible (this condition seems to be natural since
$\PGL_2(F) $ has a non-trivial $\ell$-torsion elements in this
case). 

\begin{lem}\label{lem-auxiliary1}[existence of an auxiliary place] Assume that $\ell \neq 2$ and
$\bar
\rho$ is absolutely irreducible. Assume moreover that 
$
\bar \rho \vert_{F(\zeta_\ell)}$ is absolutely irreducible if $[F(\zeta_\ell):F]= 2$. 
 Then there are infinitely many finite places $y$ such that $q_y \not \equiv 1 (\mod
\ell)$ and for the eigenvalues $\bar \alpha _y$, $\bar \beta_y$ of $\bar \rho (\Fr _y)$,
$\bar \alpha _y \neq \ q _y ^{\pm 1} \bar \beta_y $ holds.
\end{lem}
As in \cite{DT 2}, this follows from the Chebotarev density theorem and the following
lemma \ref{lem-auxiliary2}, which we give a proof since the linear disjointness of $F$ and
$\Q(\zeta_\ell)$ is assumed in the reference. The proof given here does not
use the classification of subgroups of
$\GL_2(\F_{\ell ^n})$.

 In the following, we only consider representations defined over a field of characteristic
different from $2$. We say an absolutely irreducible two dimensional representation
$\rho$ is monomial if it induced from an index
$2$ subgroup of
$G$. An absolutely irreducible representation $\rho$ is monomial if the restriction $ \rho
\vert _H$ to a normal subgroup
$H$ of
$G$ is a sum of two distinct characters. Equivalently, an
absolutely irreducible representation
$\rho$ is monomial if
$\ad ^0
\rho$ is absolutely reducible. 
\begin{lem}\label{lem-auxiliary2}
Let $k$ be a field of characteristic $\ell \neq 2$, $G$ be a finite group, 
$\rho :G
\to
\GL_2 (k)
$ be an absolutely irreducible representation, and $\chi :
G
\to k ^\times$ be an even order character. Assume that
\begin{itemize}
\item[$(\ast)$] For any $g \in G$ with
$\chi (g) \neq 1$, $\rho (g) $ has eigenvalues of the form $\{ \alpha ,\ \chi (g)
\alpha \} $.
\end{itemize}
Then $\chi$ has order $2$, and $\rho $ is induced from a character of $H$. Here $H = \ker
\chi
$.
\end{lem}
\begin{proof}[Proof of lemma \ref{lem-auxiliary2}]
We enlarge $k$ so that all eigenvalues of $\rho (g),\ g \in G, $ belong to $k$. Put $Z
=\{g \in G ,\ \rho (g)
\text{ acts as a scalar}
\}
$. By Schur's lemma, $Z$ is the center of $G$. 
$Z$ is a normal subgroup of $ H = \ker \chi $, since $\chi$ factors through $G \to G/Z
\to k ^\times$ by $(*)$. $H \neq Z$ since $\chi$ is a character of an even order. 

We need the following sublemma. 
\begin{sublem} \label{sublem-auxiliary1} Assumptions are as in \ref{lem-auxiliary2}. Then the following hold:
\begin{itemize}
\item[(a)] If $\chi (g )
$ has order $d>1$,  $g ^d \in Z$.
\item[(b)] $H/Z$ is an abelian group. 
\item[(c)] Let $ G ' = \chi ^{-1} ( \{ \pm 1\} ) $. Then $G'/ Z$ is also abelian if $ H/Z$ is
of type $(2, \ldots, 2)
$. 
\end{itemize}
\end{sublem}
\begin{proof}[Proof of \ref{sublem-auxiliary1} ] For (a), note that $\rho (g)
$ is semi-simple, and
$\rho (g )^d$ is a scalar matrix.  For (b), take an
element
$ c
\in G '$ such that $\chi (c) = -1
$. For any $h \in H$, $(c\cdot h ) ^2 \in Z$, especially $ c ^2 \in Z$, by (a). This
means that in
$H/Z$,
$c h c^{-1} = h ^{-1}$ holds, and hence $H/Z$ is abelian because the map $h \mod Z
\mapsto h ^{-1} \mod Z$ is a group homomorphism. $G'$ is a semi-direct product
of $\Z/ 2 $ by $ H/Z$. (c) is clear since the adjoint action of $c$ is trivial on $H/Z$
by the assumption.
\end{proof}

We return to the proof of lemma \ref{lem-auxiliary2}. We show $ \rho \vert _H $ is reducible. Assume contrary. Then the adjoint representation
$
\ad ^0
\rho
\vert_{H}
$ is reducible, and the irreducible consitituants are one dimensional characters since it
factors through an abelian group
$H/Z$. Thus
$\rho
\vert_H$ is monomial.
$H/Z$ is isomorphic to $\Z/ 2 $ or $\Z/ 2 \times \Z/2 $ since  it is dihedral
and abelian. Since $H/Z$ is of type $(2, \ldots, 2)
$, $G'/Z$ is also abelian by \ref{sublem-auxiliary1} (c). By the same argument, $G'/Z$ is $\Z/ 2 $ or
$\Z/ 2 \times \Z/2 $, and hence $G'/Z $ is $\Z/
2 \times \Z/2 $, and $H/Z $ is $\Z/2$. In this case, it is  easy to see that $\ad ^0
\rho \vert _{G'}
$ is a sum of three different non-trivial characters of $G'/Z$, and hence $\rho \vert _H $ can
not be irreducible, which leads to a contradiction.
\medskip

Next we show that $\rho \vert _H $ is a sum of two distinct characters. If it is
indecomposable, the unique one dimensinal subrepresentation of $H$ is fixed by the
action of
$G$, and $\rho$ is reducible. So $\rho \vert _H $ is decomposable, say $ \chi _1 \oplus 
\chi_2$. We show  $\chi_1 \neq  \chi_2$. If not,  $ \ad ^0
\rho
\vert _H 
\simeq 1  ^{
\oplus 3} 
$. 
$(*)$ implies 
$$
\tr  (g ;  \ad ^0 \rho ) = \tr ( g , 1 \oplus \chi \oplus \chi ^{-1} )
 $$
 for $
g
\in G \setminus H$, and hence  $(*)$ holds for any $ g \in G$.
$\rho
$ is reducible by the argument of \ref{prop-type1}, which leads to a contradiction.
\medskip

We finish the proof of \ref{lem-auxiliary2}. Since $\rho \vert _H  $ is a sum of two distinct
characters, $ \rho$ is monomial. $G/Z$ is dihedral, and the image of
$\chi$ is a quotient of a dihedral group.
 This implies that $\chi$ has order $2$, and  $H $ is an index $2$ subgroup of $G$.
$\rho$ is induced from a character of $H$. 
\end{proof}

\subsection{Proof of theorem B }
We prove theorem B in the introduction. 
\begin{thm} \label{thm-theoremB1}[Theorem B] Let
$ \bar \rho : G_F \to \GL _2 (\bar k_\lambda ) $ be a continuous absolutely
irreducible mod
$\ell$-representation satisfying A1)-A4):
\begin{itemize}
\item [A-1)] $ \ell \geq 3$, and $ \bar \rho \vert _{F (\zeta _{\ell} ) }
$ is absolutely irreducible if $ [F (\zeta _{\ell}) : F]= 2$.

\item [A-2)] $\bar \rho \simeq \rho _{\pi, \lambda} \mod \lambda $, $\pi $ a
cuspidal representation of $\GL_2(\A _F )$ of infinity type $ (k, w ) $ defined over
$E_\lambda$, satisfying
$ (\pi ^{\infty } ) ^K \neq \{ 0 \} $ for some compact open subgroup $ K =\prod
_u K _u $. 

\item [A-3)] for a place $ v \not \vert \ell$, $\bar \rho$
is either ramified at $v$, or
$q_v \equiv 1 \mod \ell $.

\item [A-4)]If
$g$ is even, assume that 
$\pi _u $ is essentially square integrable for some $u \neq
v $, $v \not \vert\ell$, with $ K _u = K _1 ( m ^{\cond \pi   _u}_u )$. 
\end{itemize}
Then  there is a cuspidal representation $\pi' $,  having the same
infinity type
$(k, w)$ as $\pi$, defined over a finite extension $E '_{\lambda'}  \supset E_\lambda$
such that the following conditions hold:
\begin{enumerate}
\item  The associated $\lambda ' $-adic representation $\rho_{\pi ' , \lambda '}$
gives
$\bar \rho$. 
$ (\pi ^{ ' v,
\infty} ) ^{ K ^v}\neq
\{ 0
\}
$,

\item  The conductor 
$\cond (\pi'  _v)
$ of $\pi' _v$ is equal to $ {\Art \bar \rho \vert _{G_v} }$, where $G_v$ is the
decomposition group at $v$, and $\Art
\in
\Z$ means the Artin conductor,

\item $\det
\rho _{\pi' ,
\lambda  '}\cdot  \chi ^{w +1 } _{\cycle} $ is the Teichm\"uller lift
of
$
\det
\bar
\rho \cdot  \bar \chi ^{w +1} _{\cycle}$.
\end{enumerate}

\end{thm}
\begin{proof} We prove (3) in the next subsection (proposition \ref{prop-perfect1}).  Under the assumption of
theorem B, $
\bar
\rho =
\rho _{\pi,
\lambda }
\mod
\lambda
$ for a cuspidal representation $\pi = \otimes _u
\pi _u
$ with
$ (\pi ^{ \infty} ) ^K
\neq
\{ 0 \} $ for a compact open subgroup $K = K _v \cdot K ^v$ of $\GL _2 ( \A
^\infty  _F)$. Take an integer $n \geq 0$ so that
$ K ' = K ( m_v ^n ) 
\cdot K ^v
\subset K $. $
\bar
\rho ^{I_v}
\neq
\{0\} $.  We seek for a cuspidal representation $\pi'
$ with $\mathcal O_{\lambda ' } $-coefficient 
satisfying $(\pi ^{ ' \infty  })^{K '} \neq  \{ 0 \} $,
and $\rho ^{I_v}  _{\pi' ,\lambda '}\neq \{0 \}$.

We make use of a Shimura curve. Fix $ \iota _1 \in I_F$. Choose a division quaternion
algebra
$D$ which is unramified outside $\iota_1$ when $g$ is odd, unramified outside $u$ and
$\iota _1$ when
$g$ is even. 
$ S_{\tilde K}$ denotes the canonical model over $F$ of the Shimura curve associated to
$D$ and a compact open subgroup
$\tilde K$ of
$D ^{\times } ( \A ^{\infty} ) $.
Our choice of
$\tilde K
$ is as follows.

By a Chebotarev density argument  using lemma \ref{lem-auxiliary1} and \cite {H 1} lemma 7.1, we can take an
auxiliary place
$y$ so that
$ g^{-1} K_{11} (y ) g 
\cap
\SL ( D)_+ ,\ g \in D^{\times} (\A ^\infty)$, are torsion free, $q_y \not \equiv 1 \mod
\ell$, and
$\bar 
\alpha _y
\neq q_y ^{\pm 1} \bar \beta _y
$.  Then we put

$$
\tilde K =\begin{cases}  K' \cap K_{11} (y)  \quad ( g : \text{ odd})\\
\tilde K =(\tilde K _u \cdot K ^{' u}) \cap K_{11} (y)  \quad ( g : \text{ even}).
\end{cases}
$$ 
 Here we fix a maximal order  $o_{D_u } $ of $D_u$, and 
$\tilde K _u = 1 +
\Pi ^{\cond \pi_u} _u 
\cdot o  _{D_u }$where 
$\Pi _u
$ is a uniformizer of $o _{D_u} $. 

\medskip
We take a finite set of finite places $\Sigma$ containing $\{v , v \vert \ell \}$, $\{v,   v \vert
\cond \pi \} $, $y $ and $u$ when $g$ is even.  We put $T_\Sigma = H ( D ^\times (\A
^{\Sigma, \infty }) , \tilde K )
$.

By replacing $E_\lambda$ by a bigger  $\ell$-adic field $E'_\lambda$, $
\mathcal F _{ k, w} $ defined in \ref{subsec-local}, we have the decomposition
$$
H ^1 _{\et}( S_{\tilde  K, \bar F} ,\ \bar {\mathcal F}_{ k, w}  \otimes _{\mathcal O
_\lambda} E'_{\lambda }) =
\oplus _{\pi }
\rho _{
\pi,
\lambda }\otimes _{\mathcal O _\lambda }  (\pi ^{\infty} ) ^{\tilde K }.
$$ 
Here we identify the index set as a set of cuspidal representations of $\GL_2 (\A_F)$ (not of
$D ^\times (\A) $) by the Jacquet-Langlands correspondence. By our assumption \ref{thm-theoremB1} A-2) and
 \ref{thm-theoremB1} A-4), we may assume that
$\rho  _{\pi ,
\lambda }
$ appears in
the decomposition of $ H ^1 _{\et}( S_{\tilde K, \bar F },\ \bar { \mathcal {F}} _{k, w} ) $.

By
$X$ we mean the arithmetic model of
$S_{\tilde K}
$ over
$o_v$ defined by Carayol \cite{Car 1} using Drinfeld basis, $\mathcal F = \mathcal F _{ k,w}
$.
$I = I _v $ is the inertia group at $v$.\par

Note that the arithmetic model $X$ is available for our $\tilde K$:  
In \cite{Car 1}, Carayol assumes that the compact open subgroup defining a Shimura curve is
sufficiently small. 
 Note that this condition of smallness is satisfied if we replace $ \tilde K $ by a
smaller subgroup $ U$ with $ U_v = \tilde K _v$, so it is also true for our
$\tilde K$ since $S_{\tilde K}$ is obtained by a
quotient of $S_U$ by a finite group which acts freely by our choice of auxiliary place
$y$. See also the remark on page 62 of \cite{Ja 1}.

For a finite
$\mathcal O_\lambda$-algebra
$R$, we put
$$
H _R = H ^1_{\et} ( X_{\bar
F_v}, \ \bar {\mathcal {F} } \otimes _{\mathcal O_\lambda } R) . 
$$
Note that $H_R$ carries a natural $ \mathcal O _\lambda
[ I ]
\times T_\Sigma$-module structure given by the standard action of Hecke operators.
\par 

\begin{lem} \label{lem-theoremB1}
For a maximal ideal $m$ of $T$ which is not of type $\omega$, $H _{\mathcal O_\lambda, m
}$ is
$\mathcal O_\lambda$-free, and
$$
H ^I _{\mathcal O_\lambda , m} \otimes _{\mathcal O_\lambda} k_\lambda =  H ^I _{k_\lambda , m
}
$$ holds.
\end{lem}
\begin{proof}[Proof of lemma \ref{lem-theoremB1}]
We check the conditions of lemma \ref{lem-cohomology2}, by taking $\mathcal E $ to be the category $\mathcal C_{\Omega}$ of modules of
type $\omega$.

\begin{sublem}\label{sublem-theoremB1}
The $T_\Sigma$-actions on
$ H ^0 _{\et}( X_{\bar F _v } ,\ \bar {\mathcal F} \otimes _{\mathcal O_\lambda}k_\lambda ) $ and $ H
^0 _{\et}( X_{\bar F _v } ,\
\bar {\mathcal F} ^{\vee} \otimes _{\mathcal O_\lambda} k_\lambda )$ are both of type $\omega$.
\end{sublem}
\begin{proof}[Proof of sublemma \ref{sublem-theoremB1}] The sheaf $ \mathcal F \otimes _{\mathcal O_\lambda} k_\lambda $ is trivialized in a $ D
^\times (\A ^{\Sigma , \infty } )$-equivariant way by a finite covering $X ' _{F_v} \to
X_{F_v} 
$ corresponding to $\tilde K ' = \prod _{ v\vert \ell}  K ( m_v ) \cdot \tilde K ^{\ell} 
$ by the definition \ref{subsec-local}. Over
$X'_{F_v}
$, $ D ^\times (\A ^{\Sigma, \infty} ) $-action induces the action of the
convolution algebra of type $\omega$, since it is obtained, on any constituants,
from one dimensional actions of 
$ D ^\times (\A ^{\Sigma, \infty} ) $ by 
$$\pi _ 0 ( X_{\bar F_v} ) \simeq  \pi _
0 ( F ^\times \backslash {\A ^\infty } ^\times _F/ \det \tilde K ' )
.$$
 So the claim follows for $ H ^ 0$. Since the standard action on $ H ^ 2 (X_{\bar
F_v} ,\ 
\bar
{\mathcal F}
\otimes _{\mathcal O_\lambda}k_\lambda ) $ is obtained from the dual action on $ H ^ 0 ( X_{
\bar F_v } , \ \bar { \mathcal F} _{ k, -w}\otimes _{\mathcal O_\lambda}k_\lambda  ) $ by
the Poincar\'e  duality, by the same argument the claim also follows for $ H ^ 2$. 
\end{proof}

By the proof of lemma \ref{lem-cohomology2}, we have that $H_{\mathcal O_\lambda} $ is $\mathcal O_\lambda$-free
ignoring
modules of type $\omega$. 
$ H _{\mathcal O_\lambda}\otimes _{\mathcal O_\lambda} k_\lambda = H _{k _\lambda}$ in $ \mathcal C
_{N\Omega}$ since the $T_\Sigma$-actions on
$ H ^0 _{\et}( X_{\bar F _v } ,\ \bar {\mathcal F } \otimes _{\mathcal O_\lambda}k_\lambda ) $ and $ H
^0 _{\et}( X_{\bar F _v } ,\
\bar { \mathcal F } ^{\vee} \otimes _{\mathcal O_\lambda} k_\lambda )$ are of type $\omega$.

By \cite{Car 1} 9.4.3, the arithmetic model $X= X_{ K (v ^n ) _v
\cdot K ^v }$ over
$S=
\Spec o_v$ of
$S_{K (v ^n )_v  \cdot K ^v  } $ is regular, proper and flat over $S$, if $K ^v$
is sufficiently small.
By the same holds for our $\tilde K$.

Here the set of irreducible components
$J _{ K ^v}
$ of the geometric special fiber
$X_{\bar s }$  is isomorphic to $ \oplus _{ L \in \P ^1 ( o_v / v ^n) } Y_{ L ,
K^v} $, $Y_{ L , K ^v}
$ is a smooth curve, and there is a $D ^{\times  }( \A ^{v , \infty } )$-equivariant
isomorphism
$$
\pi _0 ( Y _{L, K ^v} )\simeq
\pi _0 ( F^{\times}
\backslash
{\A_F^\infty}  ^{\times} / o ^{\times}_v \cdot \det  K ^v ).
$$ From this description, it
follows that conditions of b), c) of lemma \ref{lem-cohomology2} are satisfied for $\mathcal F  = \mathcal F_{
k,w}
$, since
$\mathcal F/
\lambda ^n \mathcal F$ is trivialized in an equivariant way by a finite covering $Y _{L,
K ^{' v}} \to Y _{L, K} $, by the same argument as above.

So lemma \ref{lem-cohomology2} is applied, and lemma \ref{lem-theoremB1} is shown. 

\end{proof}
We return to the proof of theorem \ref{thm-theoremB1}. By lemma \ref{lem-theoremB1}, $\bar \rho $ appears in $H ^I _{\mathcal O_\lambda , m} \otimes _{\mathcal
O_\lambda} k_\lambda $, and for the corresponding maximal ideal $m$, which is not of
type $\omega$, 
 $T_m $-module $H _{ \mathcal O_\lambda , m } ^I =H^I_{\mathcal O_\lambda}
\otimes _T T_m $, $H^I_{\mathcal O_\lambda}$ localized at
$m$, is
$\mathcal O_\lambda
$-free and non-zero, implying that there is a cuspidal representation $\pi ' $
with coefficient in $\mathcal O ' _{\lambda'}$ such that
$(\pi ^{ '
\infty}) ^{K'}
\neq \{ 0 \}$,
$\rho _{\pi ' ,
\lambda '}
$ gives
$\bar
\rho$ and 
$\rho ^{I_v} _{\pi ' ,
\lambda '}\neq \{  0 \} $. 

The $y$-component of $\pi '$ is spherical: since $ \pi ' _y$ has a non-zero fixed vector
by $ K_{11} (m_y )$, $\pi  '_y $ belongs to principal series or (twisted) special
representation. By our condition on $y$, the latter case does not occur since $\Fr
_y$-eigenvalues satisfies $ \bar \alpha _y \neq q_y ^{\pm 1} \bar \beta _y$. By
condition $q_y
\not
\equiv 1 \mod \ell$, $\pi ' _y$ is spherical. 

By \cite{Car 2}, th\'eor\`eme A, the Artin conductor of $\rho ' _v =
\rho _{\pi '  ,
\lambda}
\vert _{G_v}$ and the conductor of $\pi ' _v$ is the same for $ v \not \vert \ell$. By
the formula for Artin conductors,  
$$
\Art \rho'_v  = 2- \dim _{E' _{\lambda ' } } {\rho'_v} ^{I_v} + \sw \rho ' _v 
$$
$$
\Art \bar \rho \vert _{G_v } = 2- \dim _{k _{\lambda ' } } {\bar \rho } ^{I_v} + \sw
 \bar \rho \vert _{G_v }
$$
hold. Here $\sw$ means the swan conductor. Since the swan conductor does not change
under $\mod \ell$-reduction, $\Art \rho'_v  = \Art \bar \rho \vert _{G_v }$ if and only
if $ \dim _{E' _{\lambda ' } } {\rho'_v} ^{I_v} =\dim _{k _{\lambda ' } } {\bar \rho }
^{I_v} $. 
We conclude that 
$\pi'
$ satisfies the desired equality 
$$ \cond (\pi '  _ v) =\Art \bar \rho \vert
_{G_v} $$
 since $\bar
\rho$ is ramified at $v$. So $\pi ^{ ' \infty}$ has a non-zero
$ K_1 (m_v ^{\Art \bar \rho \vert _{G_v} })  \cdot \tilde K ^{ v}$-fixed vector in
the ramified cases.
\end{proof}
\bigskip
We have proved theorem B, except for the claim on the determinant.

\begin{rem}\label{rem-theoremB1}
\item{a)} The equality $\mathcal H  _R = H ^I _R$ is seen as a local invariant cycle
theorem for characteristic $\ell$ coefficient. We are mimicking the proof in the $\Q_{\ell}$-case, using
Carayol's result for
$K (v ^n ) $ that the
corresponding arithmetic model is regular, plus the determination of the ad\'ele
action on the set of irreducible components of the special fiber.
The regularity assures a purity property, in our case the Zariski-Nagata's purity
theorem for \'etale coverings suffices.

\item{b)} Even in the unramified case, we get a cuspidal representation $\pi '$
whose
$v$-component
$\pi '  _v$ has conductor at most one. By the determinant normalization in \ref{subsec-perfect}, we may replace $\pi '$ again, and obtain $\pi '$ with ${\pi ' _v }^{  K_0 (m_v)}
\neq
\{ 0
\}$. In this
case we analyze a filtration on $ \mathcal H _R$, and show the Mazur's principle in \ref{subsec-Mazurodd}. 

\item{c)} We may add $U _w$-operators for $ w \vert \ell$ (or their modification) to
$T_\Sigma$. By the same method, one can get a nearly ordinary $\pi '$ starting from
nearly ordinary $\pi$. 
\end{rem}

\subsection{Perfect complex argument}\label{subsec-perfect}
\begin{lem}\label{lem-perfect1}
Let $A$ be a noetherian local ring with maximal ideal
$m_A 
$ and the residue field $k_A$, $B$ be an $A$-algebra. Let 
$L
$ be a complex of
$B$-modules bounded below, having finitely generated cohomologies $H  ^i ( L ) $
as
$A$-modules. For a maximal ideal
$m$ of
$B$ above $m_A$, assume $ H ^i (L \otimes ^{\mathbb L} _A k _A) \otimes _B B_m  $, the
localization at
$m$, is zero for
$i
\neq 0$. Then
$ H ^0 ( L ) \otimes _B B_m $ is $A$-free. 
\end{lem}
\begin{proof} 
Since $B_m $ is $B$-flat, we may assume $B= B_m $ (finite generation assumption
is satisfied after localization since $m$ is above $m_A$). Then $H ^i (L\otimes^{\mathbb{L}} _Ak_A)$
is zero except $i = 0 $. By taking the minimal free resolution as $A$-complexes, the
claim follows. 
\end{proof}

\begin{lem}\label{lem-perfect2}
Let $\pi : X \to Y $ be an \'etale Galois covering with Galois group $G$. Let
$\mathcal F$ be a smooth $\Lambda $-sheaf on $Y$. Then $ R \Gamma (X, \pi ^* \mathcal F )$ is
a perfect complex of $\Lambda [G] $-modules, and
$$ 
 R \Gamma (X, \pi ^* \mathcal F ) \otimes^{\mathbb L} _{ \Lambda [ G] } \Lambda [G] / I
_G
\simeq
R\Gamma ( Y, \mathcal F ) 
$$
holds. Here $I_G $ is the augmentation ideal, and the map is
induced by the trace map.
\end{lem}
This is known (especially in the dual form) in any standard cohomology theory.
\begin{rem}
The above canonical morphism $  R \Gamma (X, \pi ^* \mathcal F ) \to
R\Gamma ( Y, \mathcal F ) $ in $D ^b _c ( \Lambda)$ obtained by forgetting the $G$-action is
given by the trace map.

In case of Shimura curves, which is our main application of lemma \ref{lem-perfect2}
, the above morphism is
compatible with the standard action of Hecke operators: it suffices to see the dual
morphism
$ R
\Gamma ( Y , \ \mathcal F ^\vee) \to R \Gamma ( X, \ \mathcal F ^\vee ) 
$ is compatible with dual action of $[ K g K ] $. But this is the standard action of $ [
K g^{-1} K ] $ by proposition \ref{prop-duality1}, and for the standard action, the compatibility is clear.  
\end{rem}

\medskip
We apply the lemma to adjust the determinant, which was proved by Carayol \cite{Car 3} in
case of
$\Q$ (see \cite{Ja 2} for the generalization). 

\begin{prop} \label{prop-perfect1}[determinant optimization]
Let $\pi $ be a cuspidal representation of infinite type $(k, w)$, having a non-zero
fixed vector under
$K _ 1 ( m_v ^n )  \cdot K ^v$. Then there is $\pi ' $, $\pi ^{' \infty} $ has a
non-zero fixed vector under $( K_ 1 ( m_v ^n) \cdot H )\cdot K ^v $. Here 
$$H= \{
\left(
\begin{array}{cc}
 \alpha &0 \\
 0 & \alpha\\
\end{array}
\right)
,\ \alpha \mod m _v ^n \text{ has an
}
\ell\text{-power order}\}.
$$
\end{prop}
\begin{proof} We take an auxiliary place $y$ by lemma \ref{lem-auxiliary1}, and replace $K$ by $K\cap
K_{11} (y)
$. Put
$X= S_{ K _ 1 ( v ^n ) _v \cdot K ^v}$, $Y= S_{  K _H ( v ^n ) _v \cdot K ^v}$, $\mathcal F =
\bar {\mathcal F}_{k,w}$. We view $X$ and $Y$ as complex varieties. 
$
\pi : X\to Y
$ is an
\'etale Galois covering with group $G= H / H' $, $H' $ is the inverse image of a
subgroup of
$(o _v/ m_v ^n ) ^{\times}  $ represented by units determined by $K ^v$.
We consider the complex $L$ associated to Godement's canonical resolution of $\pi ^*\mathcal
F $ on $X$, which calculates (usual, not \'etale) sheaf cohomology $H ^i ( X,\ \pi ^*\mathcal
F )$, and represents $ R \Gamma ( X, \ \pi ^* \mathcal F ) $ in the derived category. Put $A
= 
\mathcal O_\lambda[ G]$, $B = A [[T_u] ,\ [T_{u,u}] ,\ [ T_{u, u} ] ^{-1} ,\ u \not \in
\Sigma ]  
$. The complex $L$ admits a commuting action of Hecke operators, and of 
$G$ by the discussion in \ref{subsec-Hecke}, so we view $L$ as a complex of $B$-modules, sending
$[T_u]$,
$[T_{u,u}]$ to the corresponding standard Hecke action at $u$. 
 Then it follows that $ H
^i (  S_{ K _ 1 ( v ^n ) _v
\cdot K ^v} ,
\bar { \mathcal F}  _{ w, k } ) _m
$ is
$\mathcal O_\lambda [ G] $-free by the previous lemmas ($i = 0 $ or $1$ according to $g$ is
even or not), since $ H ^ 0( Y, \ \mathcal F \otimes _{\mathcal O_\lambda } k _\lambda  )$ and $
H ^ 2 ( Y, \ \mathcal F \otimes _{\mathcal O_\lambda } k _\lambda  )
$ are of type $\omega$ when $g$ is odd as in \ref{sublem-theoremB1}.  The
$G$-invariant part is non-zero by the freeness. 
\end{proof}
\section{Mazur principle}\label{sec-Mazur} In this section, we prove the following corollary A' in the introduction:

\begin{claim}\label{claim-corA}[Corollary A' (Mazur principle)]
Assumptions are as in theorem A. Then there exists $\pi ' $ as in theorem A when
$\bar
\rho $ is unramified at $v$, and $ q_v \not \equiv 1 \mod \ell$.
\end{claim}

\subsection{The odd degree case}\label{subsec-Mazurodd}
In this subsection, we discuss \ref{claim-corA} under condition A-1)-A-4) proved by Jarvis in \cite{Ja 1}. We
include this as a toy model for the Mazur principle in the even degree case in \ref{subsec-Mazureven}.  We assume that $ q
_v
\not
\equiv 1
\mod \ell$.
We may also
assume that 
$
\pi _v$ is an unramified special representation by remark \ref{rem-theoremB1} b).

We return to the general setting as in \ref{subsec-cohomology}. In addition to the assumption of lemma \ref{lem-cohomology2}, we make the following additional assumption
further:
\begin{ass}\label{ass-Mazurodd1}
\begin{enumerate}
\item
Any irreducible component $Y $ of $X_s$ is smooth.
\item 
$\mathcal F _Y$ is pure of weight $w$  for some integer $w$ independent of $Y$. 
\end{enumerate}
\end{ass}
We define and analyze a standard filtration $ W _R \subset \mathcal H  _R= H ^1 _{\et}( X_s ,\ \mathcal F _{X_s} \otimes _{\mathcal
O_\lambda } R)
$. Let $J$ be the set of irreducible component of  $X_s$, and $Z$ be the set of
singular points on $(X_s )_{\red}$.\par
We define a skyscraper sheaf $\mathcal G _R $ supported on $Z$ by
$$ 
0 \to \mathcal F \otimes _{\mathcal O_\lambda } R\to
\oplus _{Y \in J} i _Y
(\mathcal F  \vert _Y)\otimes _{\mathcal O_\lambda } R
\to
\mathcal G _R
\to 0.
$$ 
Then $ W_R \subset \mathcal H _R $ is defined by
$$
 0 \to H ^0 _{\et}(X_s ,\ \mathcal F_{X_s} \otimes _{\mathcal O_\lambda } R) \to
\oplus _{Y \in J }  H ^0 _{\et}(Y,\ \mathcal F_Y \otimes _{\mathcal O_\lambda } R) \to  H ^0_{\et}
(Z ,\
\mathcal G _R )   \to W_R \to 0 . 
$$ 
$$ 
\mathcal H _R/ W_R =
\oplus _{Y \in J } H ^1 (Y ,\ \mathcal F\vert _Y \otimes _{\mathcal O_\lambda } R)
$$
follows from the definition. In $\mathcal C_{N\Omega}$, we have $W_R = H ^0  (Z,\ \mathcal G _R)$.
Since $\mathcal G _{\mathcal O_\lambda } $ is $\mathcal O_\lambda$-smooth, it follows that the formation of $W_R$
commutes with base change after localization at maximal ideal $m$, and we conclude that
$\mathcal H _R / W_R $ has the same property.
\bigskip

We apply the formalism in \ref{subsec-cohomology} to Shimura curves $S_{D, \tilde K }$ and $\mathcal F _{k, w}$ with additional assumption \ref{ass-Mazurodd1}. As  in \ref{lem-theoremB1}, we take $\mathcal E $ in lemma \ref{lem-cohomology2} to be the category of $\omega$-type modules.

 We choose $D$ as in the proof of theorem \ref{thm-theoremB1}.  Our choice of compact open subgroup $\tilde K =  \tilde K^v \cdot \tilde K _v$ is $\tilde K ^ v = K ^v $,  and for $\tilde K _v$
$$
\tilde K _v= K(m _v ^ N ) \cdot A ,\ \quad A = \{ 
\left(
\begin{array}{cc}
\alpha & 0\\ 
0 &\alpha \\
\end{array}
\right)
,\ \alpha \in o_v ^\times,\ \alpha  \equiv 1 \mod
m_v ^ N
\} ,\ N  \geq 1 .
$$ 
The assumption \ref{ass-Mazurodd1} (1)
is satisfied for $K _v = K(m _v ^ N )$ by \cite{Car 1} 9.4.3, and the
action of $A/ A\cap K (m _v ^N ) $ is \'etale on the arithmetic model, so it is also
satisfied with our $\tilde K $. The assumption  \ref{ass-Mazurodd1} (2) on the weight of
$\mathcal F $ is satisfied with weight $w$, since over a quadratic extension of $F$
unramified at
$v$,
$\mathcal F$ is a subquotient of a (possibly higher) cohomology sheaf of an abelian scheme.
This follows from \cite{Car 1}, 2.6. See also \cite{S} for the discussion.

Let $T$ be the Hecke algebra generated by $ [T_u],\ [T_{u,u}],\
[T_{u, u}]^{-1}$ for
$u
\not
\in \Sigma$ over $\mathcal O_\lambda$,
$m$ be the maximal ideal of
$T$ corresponding to
$ {\bar
\rho } $. 

Let $P$ be the $p$-Sylow subgroup of $ K_v / K(m _v ^ N)$, where $ p$ is the residue
characteristic of $v$. It contains $K_{11} (m_v) / K (m _v ^ N )$. Since $p$ is
different from $\ell$, the operation of taking the $P$-invariants of $ \mathcal H _R $
commutes with any scalar extension. We define $U(p_v)  $-operator
for a fixed uniformizer $p_v$ at $v$ acting on $ \mathcal H _R ^ P$ in the usual way: it is defined by double coset $K
_0 ( m_v ) a(p_v) K _0 ( m_v ) \cdot K ^v$ (see section \ref{sec-notation} for the notation).
$U(p_v) 
$-operator thus defined commutes with
$G_F$-action since all group actions are defined over $F$. 

(By the vanishing of $H ^ * (P,\ - ) $,
$$ H _R ^ P = H ^ 1_{\et} (S_{K ' } ,\ \mathcal F _{ (k,w )} \otimes _{\mathcal O_{\lambda} } R )
_m,\  K ' = K ^ v
\cdot (K_{11} ( m_v ) A),
$$
and the $U(p_v)$-operator is equal to the one defined geometrically.)

By our assumption, $ \mathcal H _{E_\lambda } ^P \neq \{ 0\}$. 
Note that $(H _{k _\lambda} ^P)_m / m(H ^P _{k _\lambda})_m   $ as a Galois module is
isomorphic to $\bar \rho ^{\oplus
\alpha } 
$ for some $\alpha \geq 1 $ by the Eichler-Shimura relation \cite{Car 1}, 10.3 and by the
Boston-Lenstra-Ribet theorem \cite{BLR} using the irreducibility of $\bar
\rho$. As a $T$-module, it is isomorphic to $ (T/ m ) ^ {2\alpha }$.

 Assume the Hecke module
$T/m$ corresponding to $\bar \rho$ occurs in $ \mathcal H ^P  _{k_\lambda} / W_{k_\lambda} ^P
$. Then it occurs in
$ (\mathcal H ^P _{\mathcal O_\lambda } / W^P_{\mathcal O_\lambda} )_m \neq \{0 \}$, which is $\mathcal
O_\lambda$-free. Let
$\pi
$ be the corresponding cuspidal representation.  Since
$\mathcal F _Y$ is pure of weight $w$ for each irreducible component $ Y $ of $X_s$, the weight
here is
$w +1$ by the Weil conjecture. If we assume that
$\pi _v$ is special, $\det \rho _{\pi, \lambda} (\Fr _v) = q _v \beta  ^2$, with
$\vert \beta \vert = q^{\frac {w +1} 2 } _v $, which is impossible by the weight
reason. 
It follows that $\pi _v$ must be a principal series.
Moreover, $\rho _v = \rho_{\pi , \lambda
' }
\vert _{G_v} 
$, associated to $\pi_v$ by the local Langlands correspondence, has a non-zero $I_v
$-fixed part and
$\det
\rho _v
\vert_{I_v} $ is the Teichm\"uller lift of $ \det \bar \rho \vert_{I_v}$ by our choice of
$ \tilde K _v$. We conclude that $\pi_v$ is an unramified principal series. 
\par 

So we may assume all $T/m $ appear in $W_{k_\lambda}  ^P $ and hence $\bar \rho ^{\oplus
\alpha} 
\simeq W_{k_\lambda}^P/ m W_{k_\lambda}^P $ as a Galois-Hecke bimodule. For any lift
$\pi$ found in $ W_{E_\lambda}^P$, 
$\pi $ has an unramified special component at $v$. $q_v U (p_v) \cdot \Fr ^{-1} _v$ is
identity on $W^P _{E _\lambda }$ by the compatibility of local and global Langlands
correspondence recalled in section \ref{sec-notation}  (by our normalization, the eigenvalue of
$U(p_v)$-operator is equal to the Frobenius eigenvalue on $ \rho_{\pi, \lambda}  / \rho
^{I_v}_{\pi,
\lambda}
$), and hence on
$ W^P _{\mathcal O  _\lambda }$ also by the
$\mathcal O_\lambda$-freeness. Since
$U(p_v)$ commutes with global Galois action on
$\bar
\rho$, by Schur's lemma we may assume that there is a scalar $ \gamma \in \bar k
_\lambda
\setminus \{ 0\}
$ such that
$
\gamma ^{-1}  \Fr _v 
$ acts trivially on 
$$  W_{k_\lambda}^P/ m W_{k_\lambda}^P =(H _{k_\lambda} ) ^P _m / m (H
_{k_\lambda} ) ^P _m\otimes _{k_\lambda} \bar k _ {\lambda} \simeq \bar \rho ^{\oplus
\alpha}. 
$$It follows that any
$\Fr _v$-eigenvalues of $\bar \rho$ are the same, which leads to a contradiction since two
eigenvalues of $\Fr _v $ on
${\bar
\rho } 
$ are of the form $\bar \alpha _v ,\ q_v \bar \alpha _v$, and
$q_v \not \equiv 1 \mod \ell $.

\subsection{Cerednik-Drinfeld type theorem for totally real fields}\label{subsec-CD}

Let $D$ be a quaternion algebra over $F$. Assume $D$ defines a Shimura curve,
i.e., $D \otimes_{\Q} \R \simeq M _2 (\R)  \times \H ^{g-1}$.  

Choose a finite place $v$ where $\inv _ v D = 1/2$. Let
$\tilde D $ be the definite quaternion with $\inv _v D= 0$, and
other invariants at finite places are the same as $D$.  
We need a generalization of the Cerednik-Drinfeld theorem for Shimura curves over $\Q$. For general totally real fields, such a result
follows from the work of Boutot-Zink \cite{BZ}, which we recall in the following.

We choose a compact open subgroup $ K \subset D^\times(\A ^\infty ) $ such that $
K = o ^\times _{D_v} 
 \cdot  K ^v $.
Put  $\tilde K = \GL _2 (o _v ) \cdot K ^v$.
Then the main result of \cite{BZ} claims:
\begin{thm}\label{thm-CD1}[\cite{BZ}, theorem 0.1]
There is a canonical isomorphism
$$
S _{D, K,\ o _v}\simeq \tilde D ^{\times} \backslash D^\times (\A ^{v,
\infty})\times
\hat
\Omega _{o _v}\hat \otimes _{o _v }o^{\unr} _v / K ^v .
$$
Here $\hat \Omega _{o
_v}$ is the Deligne model of the Drinfeld upper half plane, and $o^{\unr} _v $ the
maximal unramified extension of $o_v$. The action of
$\GL _2 (F_v) $ on $ \hat \Omega _{o
_v} \hat \otimes _{o _v }o^{\unr} _v $ is 
$$
g \mapsto ( g ,\ \Fr ^{n(g)} _v  ),
 \quad n (g) = \ord _v (\det \ g) , \ \Fr _v : 
\text{ Frobenius at } v
$$
and $\tilde D ^{\times}$ acts diagonally. Especially, $S _{D, K,\ o _v}$ is regular.
\end{thm}

\begin{cor}\label{cor-CD1}
The set of irreducible components $J _K $of $ S _{K,\ o ^{\unr}  _v}$ is identified
with two copies of the double cosets
$S _{\tilde D, \tilde K }= \tilde D ^{\times}
\backslash \tilde D ^{\times} ( \A ^{\infty} ) / \tilde K 
$, and the identification is $\tilde D ^\times  (\A^{v, \infty} )$-equivariant.
\end{cor}
\begin{proof} The set of irreducible components $I$ of $\hat \Omega _{o_v}  $ is
isomorphic to the set of all lattices up to homothety in $ F _v ^2 $ by the structure of
the Bruhat-Tits building, and hence identified with
$
\GL_2 (F_v) /
F_v ^\times \cdot \GL_2 (o_v) 
$. 

By theorem \ref{thm-CD1}, $J_K $ is canonically isomorphic to 
$$\tilde D ^{\times} \backslash  D^\times (\A ^{v, \infty}) \times I \times  \Z / K ^v
.
$$  
Put $ \tilde D ^\times _+ = \{g \in \tilde D ^\times  ,\ \ord (\det g) \text{ is even} \}
$.

The map $\alpha :  \tilde D ^\times _+ \to 2 \Z  $ given by $ \ord (\det g) $ is
surjective. So the pieces
$J_{K, + } = \tilde D ^{\times}_+ 
\backslash  D^\times (\A ^{v, \infty}) \times I
\times  2\Z / K ^v$ and $J_{K , -} = 
\tilde D ^{\times}_+  \backslash  D^\times (\A ^{v, \infty}) \times I \times  ( \Z
\setminus 2\Z )  / K ^v$ inside $J_K$ are isomorphic to $  \alpha ^{-1}
(0) \backslash\GL_2 (F_v) / (F_v ^\times \cdot \GL_2 (o_v) )/ K ^v \simeq S_{\tilde D,
\tilde K } $.
$\tilde D ^\times  (\A^{v,
\infty} )$-equivariance follows from the description. Moreover, the construction of $J
_{K \pm } $ is canonical.
\end{proof}
\medskip
We need to calculate the fibers of our sheaves $\mathcal F^D _{k,w}$ as well. Note that we
also have a sheaf $ \mathcal F ^{\tilde D}_{k,w} $ which is an analogue of $\mathcal F ^D
_{k,w}$ on
$S _{D, K}$.
\begin{cor}\label{cor-CD2}
The sum of geometric generic fibers $ \oplus _{\eta \in J _K }  (\bar {\mathcal F} _{k, w
})_{\bar
\eta}$ is identified with the sheaf
$\mathcal F ^{\tilde D } _{k,w}$ on 
$S _{\tilde D , \tilde K }
$ on each $ J _{K, \pm}$, and the identification is $\tilde D ^\times  (\A^{v,
\ell, \infty} )$-equivariant.
\end{cor}
\begin{proof}
Let $\pi _{\ell} : \tilde S _{\ell} \to S _{D, K} $ be the Galois covering
corresponding to $\prod _{u \vert \ell} K _u / K \cap \overline {(F ^{\times})} $. By
corollary \ref{cor-CD1}, the set of irreducible components of $(\tilde S _{\ell}
)_{\overline{k(v)}}$ is identified with two copies of $S_{\tilde D , \tilde K ^{\ell}} $,
preserving
$ \tilde D ^{\times}(\A ^{v, \infty}) $-action. Since the sheaf $ \mathcal F ^D _{k, w}$
is obtained from the covering $\pi _{\ell}$ and the representation $ \otimes
_{\iota \in I _F} (\iota \det) ^{k' _{\iota}}\Sym ^{ k _{\iota}-2}$ of $ D ^{\times}
(\bar {\Q} _{\ell})$ by contracted product, their geometric generic fibers are identified
with $\mathcal F ^{\tilde D} _{k, w} $ on each $ J _{K , + } $ and $ J _{K , -} $ by the $
\tilde D ^{\times}(\A ^{v,
\infty})
$-equivariance.
\end{proof}

We need $U_v$-operator in the following, which is defined as follows. Take a uniformizer
$\Pi _v$ of $\mathcal O_{D_v} $, and consider the double coset $K ^ v \cdot K_v \Pi _v K_v $.
This defines a correspondence, and the action on cohomology is $U_v$. $U_v$-operator
thus defined exists over $F$, not only over $F_v$, and by the local Jacquet-Langlands
correspondence \cite{JL}, it corresponds to $U(p_v )$-operator defined in \ref{subsec-Mazurodd} for
$\GL_2 (F_v )$ (this is proved easily).

\subsection{The even degree case}\label{subsec-Mazureven}
By \ref{thm-CD1} and the method of section \ref{subsec-Mazurodd}, we prove \ref{claim-corA} without assuming A-4). 
By subsection \ref{subsec-Mazurodd} and remark \ref{rem-theoremB1} b), we may assume that the degree $[F:\Q ]$ is even, and $\pi_v$ is a special representation twisted by an unramified character. 

 Take
$\tilde D$, the definite quaternion algebra which is unramified at all
finite places, and let
$D$ be an indefinite quaternion algebra corresponding to $\tilde D$ as above.
By the argument of
\ref{subsec-auxiliary}, we take an auxiliary place $y$, and take a compact open subgroup $ K
= 
\prod _u  K _u $ with $ K _ v = o ^\times _{D_v} $, $ K _ y= K _{11} ( m
_y) $.
Then $\rho _{\pi , \lambda} $ occurs in $H ^1_{\et}(S_{D, K } ,
\bar {\mathcal F} ^D _{k,w})$ by the Jacquet-Langlands
correspondence, since $\pi _v$ is an unramified special representation.

We choose $\Sigma $ so that $u \not \in \Sigma \Rightarrow u \not \vert \ell$ and $ K_u =
\GL _2 (o _u ) $. Let $T$ be the Hecke algebra generated by $ [T_u],\ [T_{u,u}],\
[T_{u, u}]^{-1}$ for
$u
\not
\in \Sigma$ over $\mathcal O_\lambda$,
$m$ be the maximal ideal of
$T$ corresponding to
$ {\bar
\rho } $. Assume contrary, so 
$ {\bar \rho } $ does not appear in $ H ^0 ( S _{\tilde D, \tilde K}, \ \bar { \mathcal F}
^{\tilde D}  _{k, w} \otimes _{\mathcal O_\lambda} k_\lambda)
$. 
\begin{prop}\label{prop-Mazureven1}
$$H ^1 _{\et}( S_{D ,  K , o ^{\unr} _v} ,  \mathcal F
^{D}  _{k, w} \otimes _{\mathcal O_\lambda}R ) _m \simeq H ^1 _{\et}( S_{D ,
 K , \bar F_v } , \bar { \mathcal F}
^{ D}  _{k, w} \otimes _{\mathcal O_\lambda} R )^{I_v} _m
$$
holds for a finite local $\mathcal O_\lambda$-algebra $R$, and hence $R \mapsto H ^1 _{\et}(
S_{D ,
 K } , \mathcal F
^{ D}  _{k, w} \otimes _{\mathcal O_\lambda} R )^{I_v} _m$ commutes with scalar extensions.
\end{prop}
\begin{proof} We adopt the method used in previous sections. Since our $T$-action
annihilates $( \oplus _{Y \in J} H ^0_{\et}(Y, \mathcal F \vert _Y) ) _m \simeq H ^0 ( S _{\tilde
D,
\tilde K},
\  \mathcal F ^{\tilde D}  _{k, w} \otimes _{\mathcal O_\lambda} k_\lambda) ^{\oplus 2} _m $ by
our assumption that $\bar \rho$ does not come from $ H ^0 ( S _{\tilde D, \tilde K}, \ \bar { \mathcal F}
^{\tilde D}  _{k, w} \otimes _{\mathcal O_\lambda} k_\lambda)
$, the claim follows from \ref{cor-CD2} and \ref{lem-cohomology2}.
\end{proof}
The rest of the argument is the same as in section \ref{subsec-Mazurodd}.
By
proposition \ref{prop-Mazureven1}, ${\bar \rho } $ comes from $ H ^1_{\et} ( S_{D ,  K , \bar  F _v} , \ 
\bar { \mathcal F} ^{\tilde D}  _{k, w} ) ^{I_v}
_m$ as a Hecke module. As in section \ref{subsec-Mazurodd}, $q_v U_v  \cdot \Fr _v ^{-1}  $ is identity on
$ H ^1_{\et} ( S_{D ,  K , \bar F _v} , \  \bar { \mathcal F}
^{\tilde D}  _{k, w} ) ^{I_v}$, and $U_v$-operator on $ H ^1 _{\et}( S_{D ,  K , \bar F
_v} ,
\ 
\bar{ \mathcal F} ^{\tilde D}  _{k, w} \otimes _{\mathcal O _\lambda} k_{\lambda}) ^{I_v}_m \otimes
_T T/ m $ acts as a scalar since it commutes with the global Galois
group $G_F$ by using Boston-Lenstra-Ribet theorem\cite{BLR}. Two eigenvalues of $\Fr _v
$ on
${\bar
\rho } 
$ are of the form $\bar \alpha _v ,\ q_v \bar \alpha _v$, this implies $ q _v \equiv 1
\mod
\ell$.
\section{Proof of theorem A}\label{sec-theoremA}
Now we finish the proof of theorem A in the introduction. We may assume that the degree $[F: \Q]$ is
even. By 
\cite{T}, theorem 1, there is a finite place $z\neq v  $ of $F$ where $\bar \rho $ is
unramified, $ q _ z \equiv -1 \mod \ell $, and there is a cuspidal representation $\tilde \pi $
which gives $\bar \rho$ such that $\tilde \pi _z $ is an unramified special representation at $z$,
$(\tilde
\pi ^
\infty ) ^{ K
\cap K_0  (z)} \neq \{ 0 \}$. We can now apply theorem B to $\tilde \pi $, and
optimize the level at $v$. When $\bar \rho$ is unramified at $v$, we
first apply remark \ref{rem-theoremB1} b) to get $\pi '$ which has an unramified special component $\pi '_v$ at $v$, then
apply the Mazur principle in the form of section \ref{subsec-Mazurodd}. Finally, the auxiliary place
$z$ can be removed by the Mazur principle in \ref{subsec-Mazureven}, since $ q _z \not \equiv 1 \mod \ell
$, and the component at $z$ is an unramified special representation.

e-mail address: fujiwara@math.nagoya-u.ac.jp
\enddocument